\documentclass[11pt]{article}
\usepackage{amssymb,amsmath,amsthm, amsfonts}

\def\sqr#1#2{{\vcenter{\vbox{\hrule height.#2pt
              \hbox{\vrule width.#2pt height#1pt \kern#1pt \vrule width.#2pt}
          \hrule height.#2pt}}}}
\def\signed #1{{\unskip\nobreak\hfil\penalty50
          \hskip2em\hbox{}\nobreak\hfil#1
          \parfillskip=0pt \finalhyphendemerits=0 \par}}
\def\endpf{\signed {$\sqr69$}}
\def\sqr#1#2{{\vcenter{\vbox{\hrule height.#2pt
              \hbox{\vrule width.#2pt height#1pt \kern#1pt \vrule width.#2pt}
              \hrule height.#2pt}}}}
\def\signed #1{{\unskip\nobreak\hfil\penalty50
              \hskip2em\hbox{}\nobreak\hfil#1
              \parfillskip=0pt \finalhyphendemerits=0 \par}}
\def\endpf{\signed {$\sqr69$}}
\def\3n{\negthinspace \negthinspace \negthinspace }
\def\2n{\negthinspace \negthinspace }
\def\1n{\negthinspace }

\def\={\buildrel \triangle \over =}

\def\O{\Omega}

\def\q{\quad}

\def\max{\mathop{\rm max}}

\def\exp{\mathop{\rm exp}}
\def\sup{\mathop{\rm sup}}

\def\|{\Big |}
\def\({\Big (}
\def\){\Big )}
\def\[{\Big[}
\def\]{\Big]}
\def\be{\begin{equation}}
\def\bel{\begin{equation}\label}
\def\ee{\end{equation}}
\def\bt{\begin{theorem}}
\def\bcd{\begin{condition}}
\def\ecd{\end{condition}}
\def\et{\end{theorem}}
\def\bc{\begin{corollary}}
\def\ec{\end{corollary}}
\def\bde{\begin{definition}}
\def\ede{\end{definition}}
\def\bl{\begin{lemma}}
\def\el{\end{lemma}}
\def\bp{\begin{proposition}}
\def\ep{\end{proposition}}
\def\br{\begin{remark}}
\def\er{\end{remark}}
\def\ba{\begin{array}}
\def\ea{\end{array}}
\def\ed{\end{document}}

\def\square#1{\vbox{\hrule\hbox{\vrule height#1%
     \kern#1\vrule}\hrule}}
\def\rectangle#1#2{\vbox{\hrule\hbox{\vrule height#1%
     \kern#2\vrule}\hrule}}

\font\tenbb=msbm10 \font\sevenbb=msbm7 \font\fivebb=msbm5

\newfam\bbfam
\scriptscriptfont\bbfam=\fivebb \textfont\bbfam=\tenbb
\scriptfont\bbfam=\sevenbb

\newtheorem{lemma}{Lemma}[section]
\newtheorem{remark}{Remark}[section]

\newtheorem{theorem}{Theorem}[section]
\newtheorem{corollary}{Corollary}[section]

\newtheorem{definition}{Definition}[section]
\newtheorem{proposition}{Proposition}[section]
\newtheorem{condition}{Condition}[section]

\makeatletter
   
   \@addtoreset{equation}{section}
\makeatother

\begin{document}
\title{ Stochastic Differential Games and Viscosity Solutions of
Hamilton-Jacobi-Bellman-Isaacs Equations \bf\footnote{The work of
Mme LI has been supported by a one-year fellowship awarded by the
General Council of Finist$\grave{e}$re, France, the NSF of P.R.China
(No. 10426022; 10371067) and Program for Changjiang Scholars and
Innovative Research Team in University (PCSIRT).}}

\author{ Rainer Buckdahn\\
{\small D$\acute{e}$partement de Math$\acute{e}$matiques,
Universit$\acute{e}$ de Bretagne
Occidentale,}\\
 {\small 6, avenue Victor-le-Gorgeu, B.P. 809, 29285 Brest
cedex, France.}\\
{\small{\it E-mail: Rainer.Buckdahn@univ-brest.fr.}}\\
 Juan Li\\
{\small School of Mathematical
Sciences, Fudan University, Shanghai 200433,}\\
{\small Department of Mathematics, Shandong University at Weihai, Weihai 264200, P. R. China.}\\
{\small {\it E-mail: juanli@sdu.edu.cn.}}  }
\date{}
\maketitle \noindent{\bf Abstract}\hskip4mm
  In this paper we study zero-sum two-player stochastic differential
  games with the help of theory of Backward Stochastic
Differential Equations (BSDEs). At the one hand we generalize the
results of the pioneer work of Fleming and Souganidis~\cite{FS1} by
considering cost functionals defined by controlled BSDEs and by
allowing the admissible control processes to depend on events
occurring before the beginning of the game (which implies that the
cost functionals become random variables), on the other hand the
application of BSDE methods, in particular that of the notion of
stochastic ``backward semigroups" introduced by Peng~\cite{Pe1}
allows to prove a dynamic programming principle for the upper and
the lower value functions of the game in a straight-forward way,
without passing by additional approximations. The upper and the
lower value functions are proved to be the unique viscosity
solutions of the upper and the lower Hamilton-Jacobi-Bellman-Isaacs
equations, respectively. For this Peng's BSDE method
(Peng~\cite{Pe1}) is translated from the framework of stochastic
control theory into that of stochastic
differential games.
\newpage
\section{\large{Introduction}}

\hskip1cm With their pioneer paper of 1989 Fleming and Souganidis
[8] were the first to study in a rigorous manner two-player
zero-sum stochastic differential games and to prove that the lower
and the upper value functions of such games satisfy the dynamic
programming principle, that they are the unique viscosity
solutions of the associated Bellman-Isaacs equations and coincide
under the Isaacs condition. Their work has translated former
results by Evans and Souganidis [7] from a deterministic into the
stochastic framework and has given an important impulse for the
research in the theory of stochastic differential games. And so a
lot of recent works are based on the ideas developed in [8], see,
for instance, Buckdahn, Cardaliaguet, Rainer [4], Hou, Tang [11]
and Rainer [16]. The reader interested in this subject is also
referred to the references given in [8].

Also the present work investigates two-player zero-sum stochastic
differential games, but with two main differences to the setting
chosen by Fleming and Souganidis [8] and the other papers mentioned
above: At the one hand we allow our admissible control processes to
depend on the full past of the trajectories of the driving Brownian
motion, this means, in particular they can also depend on
information occurring before the beginning of the game (which has
the consequence that the cost functionals become random variables),
on the other hand we consider a more general running cost
functional, which implies that the cost functionals will be given by
a backward stochastic differential equation (for short, BSDE). These
both extensions of the framework in [8] are crucial because they
allow to harmonize the setting for stochastic differential games
with that for the stochastic control theory and to simplify
considerably the approach in [8] by using BSDE methods.

BSDEs in their general non-linear form were introduced by Pardoux
and Peng [12] in 1990. They have been studied since then by a lot
of authors and have found various applications, namely in
stochastic control, finance and the second order PDE theory. BSDE
methods, originally developed by Peng [14], [15] for the
stochastic control theory, have been introduced in the theory of
stochastic differential games by Hamad\`{e}ne, Lepeltier [9] and
Hamad\`{e}ne, Lepeltier and Peng [10] to study games with a
dynamics whose diffusion coefficient is strictly elliptic and
doesn't depend on the controls. In our present work there isn't
any such restriction on the diffusion coefficient and the
application of BSDE methods, in particular the notion of
stochastic backward semigroups (Peng [14]), allows to prove the
dynamic programming principle for the upper and lower value
functions of the game in a very straight-forward way (i.e., in
particular without making use of $r$-strategies and
$\pi$-admissible strategies playing an essential role in [8]) and
to derive from it with the help of Peng's method (see [14], [15])
the associated Bellman-Isaacs equations.

The dynamics of the stochastic differential game we investigate is
given by the controlled stochastic differential equation \be
  \left \{
  \begin{array}{llll}
  dX^{t,x ;u, v}_s & = & b(s,X^{t,x; u,v}_s, u_s, v_s) ds + \sigma(s,X^{t,x; u,v}_s, u_s, v_s) dB_s,  \\
   X^{t,x ;u, v}_t  & = & x (\in \mathbb{R}^n),\hskip 7cm s\in
   [t,T],
   \end{array}
   \right.
  \ee
where $T>0$ is an arbitrarily fixed finite time horizon,
$B=\left(B_s\right)_{s\in [0,T]}$ is a $d$-dimensional standard
Brownian motion, and $u=\left(u_s\right)_{s\in[t,T]}$,
$v=\left(v_s\right)_{s\in[t,T]}$ are progressively measurable with
respect to the Brownian filtration and take their values in some
compact metric spaces $U$ and $V$, respectively (we will say that
$u\in{\cal U}_{t,T},v\in{\cal V}_{t,T}$). Precise assumptions on
the coefficients $b:[0,T]\times \mathbb{R}^n\times U\times
V\rightarrow \mathbb{R}^n$ and $\sigma:[0,T]\times
\mathbb{R}^n\times U\times V\rightarrow \mathbb{R}^{n\times d}$
are given in the next section.

The cost functional (interpreted as a payoff for Player I and as a
cost for Player II) is introduced by a backward stochastic
differential equation (BSDE, for short):
 \be
   \left \{\begin{array}{rcl}
   -dY^{t,x; u, v}_s & = & f(s,X^{t,x; u, v}_s, Y^{t,x; u, v}_s, Z^{t,x; u, v}_s,
                              u_s, v_s) ds -Z^{t,x; u, v}_s dB_s,\\
      Y^{t,x; u, v}_T  & = & \Phi (X^{t,x; u, v}_T),\qquad
\qquad\qquad \qquad \qquad \qquad \quad s\in [t,T],
   \end{array}\right.
   \ee
where the driver $f:[0,T]\times \mathbb{R}^n\times
\mathbb{R}\times \mathbb{R}^d\times U\times V\rightarrow
\mathbb{R}$ describes the running cost and
$\Phi:\mathbb{R}^n\rightarrow \mathbb{R}$ the terminal cost. Under
the assumptions on $f$ and $\Phi$ that will be introduced in the
next section the above BSDE has a unique solution $(Y_s^{t,x;
u,v},Z_s^{t,x; u,v})_{s\in[t,T]}$ and the cost functional is given
by
\begin{equation}
J(t,x;u,v)=Y_t^{t,x; u,v}.
\end{equation}
As usual in the differential game theory, the players cannot
restrict to play only control processes, one player has to fix a
strategy while the other player chooses the best answer to this
strategy in form of a control process. A strategy admissible for
Player I (resp., Player II) is a non-anticipating mapping
$\alpha:{\cal V}_{t,T}\rightarrow {\cal U}_{t,T}$\ (resp.,
$\beta:{\cal U}_{t,T}\rightarrow {\cal V}_{t,T}$) which associates
every admissible control of the other player with one of his own
admissible controls (we write: $\alpha\in{\cal
A}_{t,T},\beta\in{\cal B}_{t,T}$; the precise definitions can be
found in Section 4). We define the lower value function of our
stochastic differential game as follows:
\begin{equation}
    W(t,x):=\text{essinf}_{\beta\in{\cal B}_{t,T}}
    \text{esssup}_{u\in{\cal U}_{t,T}}
    J(t,x;u,\beta(u))
\end{equation}
and the upper value function is given by
\begin{equation}
     U(t,x):=\text{esssup}_{\alpha\in{\cal A}_{t,T}}
    \text{essinf}_{v\in{\cal V}_{t,T}}
    J(t,x;\alpha(v),v).
\end{equation}
\smallskip

The objective of our paper is to investigate these lower and upper
value functions.  The main results of the paper state that $W$ and
$U$ are deterministic (Proposition 4.1) continuous viscosity
solutions of the Bellman-Isaacs equations (Theorem 5.1)
 \be
 \left \{\begin{array}{ll}
 &\!\!\!\!\! \frac{\partial}{\partial t}  W(t,x) +  H^{-}(t, x, W, DW, D^2W)=0,
 \hskip 0.5cm (t,x)\in [0,T)\times {\mathbb{R}}^n ,  \\
 &\!\!\!\!\!  W(T,x) =\Phi (x), \hskip0.5cm x \in {\mathbb{R}}^n,
 \end{array}\right.
\ee
and
\be
 \left \{\begin{array}{ll}
 &\!\!\!\!\! \frac{\partial }{\partial t} U(t,x) +  H^{+}(t, x, U, DU, D^2U)=0,
 \hskip 0.5cm (t,x)\in [0,T)\times {\mathbb{R}}^n ,  \\
 &\!\!\!\!\!  U(T,x) =\Phi (x), \hskip0.5cm  x \in {\mathbb{R}}^n,
 \end{array}\right.
\ee
respectively, associated with the Hamiltonians
$$ H^-(t, x, y, p, X)=
\mbox{sup}_{u \in U}\mbox{inf}_{v \in V}H(t,x,y,p,X,u,v),$$
$$ H^+(t, x, y, p, X)= \mbox{inf}_{v \in
V}\mbox{sup}_{u \in U}H(t,x,y,p,X,u,v),$$ $(t,x,y,p,X) \in
[0,T]\times \mathbb{R}^n\times \mathbb{R}\times \mathbb{R}^n\times
S^n$ (Recall that $S^n$ denotes the set of all $n\times n$
symmetric matrices), where
\begin{eqnarray}
\nonumber H(t,x,y,p,X,u,v) &=&
1/2\cdot\text{tr}\left(\sigma\sigma^T(t,x,u,v)X \right)\\
&+& p\cdot b(t,x,u,v)+f(t,x,y,p\cdot \sigma(t,x,u,v),u,v).
\end{eqnarray}
Moreover, we prove the uniqueness (Theorem 6.1) in a class of
continuous functions with a growth condition which was introduced
by Barles, Buckdahn and Pardoux [3] and is weaker than the
polynomial growth assumption.

Notice that the fact that $W$ and $U$, introduced as combination of
essential infimum and essential supremum over a class of random
variables, are deterministic is far from beng trivial. The method
developed by Peng [14, 15] (see also Theorem 3.1 of the present
paper) for value functions involving only control processes but not
strategies doesn't apply here since the strategies from ${\cal
A}_{t,T}$\ and ${\cal B}_{t,T}$ don't have, in general, any
continuity property. To overcome this difficulty we show in
Proposition 4.1 and Lemma 4.1 that $W$ and $U$ are invariant under
Girsanov transformation and use the fact that a functional of the
Brownian motion which is invariant under Girsanov transformation
into all directions of the Cameron-Martin space must be
deterministic. We emphasize that the proof of Lemma 4.1 doesn't use
BSDE methods which makes this method also applicable to the other
situations, such as standard stochastic control problems.

Our paper is organized as follows. The Sections 2 and 3 recall
some elements of the theory of backward SDEs and forward-backward
SDEs which will be needed in the sequel. Section 4 introduces the
setting of the stochastic differential game and its lower and
upper value functions $W$ and $U$, and proves that these both
functions are deterministic and satisfy the dynamic programming
principle (for short, DPP). The DPP allows to derive in Section 5
with the help of Peng's method that $W$ and $U$ are viscosity
solutions of the associated Bellman-Isaacs equations; the
uniqueness is studied in Section 6. Finally, after having
characterized $W$ and $U$ as unique viscosity solutions of
associated Bellman-Isaacs equations we show that under the Isaacs
condition $W$ and $U$ coincide (one says that the game has a
value) and we also identify $W$ and $U$ with the value functions
defined in [8].

\section{ {\large Preliminaries}}

  \hskip1cm Let us begin by introducing the setting for the stochastic
differential game we want to investigate. We consider as Brownian
motion B is the d-dimensional coordinate process on the classical
Wiener space $(\Omega, {\cal{F}}, P)$, i.e., $\Omega$ is the set of
continuous functions from [0, T] to ${\mathbb{R}}^d$ starting from 0
($\Omega= C_0([0, T];{\mathbb{R}}^d)$), $ {\cal{F}} $ the completed
Borel $\sigma$-algebra over $\Omega$, P the Wiener measure and B the
canonical process: $B_s(\omega)=\omega_s,\ s\in [0, T],\ \omega\in
\Omega$. By $\{{\mathcal{F}}_s,\ 0\leq s \leq T\}$\ we denote the
natural filtration generated by $\{B_s\}_{0\leq s\leq T}$\ and
augmented by all P-null sets, i.e.,
$${\mathcal{F}}_s=\sigma\{B_r, r\leq s\}\vee {\mathcal{N}}_P,\ \  s\in [0, T], $$
where $ {\cal{N}}_P$ is the set of all P-null subsets, and $T > 0$\
a fixed real time horizon. For any
    $n\geq 1,$\ $|z|$ denotes the Euclidean norm of $z\in
    {\mathbb{R}}^{n}$. We also shall introduce the following both spaces of processes which will be used frequently in the sequel:
\vskip0.2cm
    ${\cal{S}}^2(0, T; {\mathbb{R}}):=\{(\psi_t)_{0\leq t\leq T}\mbox{ real-valued adapted c\`{a}dl\`{a}g
    process}:\\ \mbox{ }\hskip6cm
    E[\sup\limits_{0\leq t\leq T}| \psi_{t} |^2]< +\infty \}; $
    \vskip0.2cm

   ${\cal{H}}^{2}(0,T;{\mathbb{R}}^{n}):=\{(\psi_t)_{0\leq t\leq T}\ {\mathbb{R}}^{n}\mbox{-valued progressively
   measurable process}:\\ \mbox{ }\hskip6cm
     \parallel\psi\parallel^2_{2}=E[\int^T_0| \psi_t| ^2dt]<+\infty \}. $
\vskip0.2cm
 Let us now consider a function $g:
\Omega\times[0,T]\times {\mathbb{R}} \times {\mathbb{R}}^{d}
\rightarrow {\mathbb{R}} $ with the property that $(g(t, y,
z))_{t\in [0, T]}$ is progressively measurable for each $(y,z)$ in
${\mathbb{R}} \times {\mathbb{R}}^{d}$, and we also make the
following assumptions on $g $ throughout the paper:
 \vskip0.2cm

(A1) There exists a constant $C\ge 0$  such that, P-a.s., for all
$t\in [0, T],\ y_{1}, y_{2}\in {\mathbb{R}},\ z_{1}, z_{2}\in
{\mathbb{R}}^d,\\ \mbox{ }\hskip4cm   |g(t, y_{1}, z_{1}) - g(t,
y_{2}, z_{2})|\leq C(|y_{1}-y_{2}| + |z_{1}-z_{2}|).$
 \vskip0.2cm

(A2) $g(\cdot,0,0)\in {\cal{H}}^{2}(0,T;{\mathbb{R}})$.
\vskip0.2cm

 The following result on backward stochastic differential equations (BSDEs) is by now well known, for its proof the reader is referred to
 Pardoux and Peng~\cite{PaPe}.
 \bl Under the assumptions (A1) and (A2), for any random variable $\xi\in L^2(\O, {\cal{F}}_T,$ $P),$ the
BSDE
 \be y_t = \xi + \int_t^Tg(s,y_s,z_s)ds - \int^T_tz_s\,
dB_s,\q 0\le t\le T, \label{BSDE} \ee
 has a unique adapted solution
$$(y^{T, g, \xi}_t, z^{T, g,
\xi}_t)_{t\in [0, T]}\in {\cal{S}}^2(0, T; {\mathbb{R}})\times
{\cal{H}}^{2}(0,T;{\mathbb{R}}^{d}). $$ \el
 In the sequel, we
always assume that the driving coefficient $g$\ of a BSDE
satisfies (A1) and (A2).

Let us remark that Lemma 2.1 remains true when assumption (A1) is
replaced by weaker assumptions, for instance those studied in
Bahlali~\cite{B}, Bahlali, Essaky, Hassani and Pardoux~\cite{BEHP}
or Pardoux and Peng~\cite{PP}. However, here, for the sake of
simplicity of the calculus we prefer to work with the Lipschitz
assumption.

   We also shall recall the following both basic results on BSDEs.
   We begin with the well-known comparison theorem (see
   El Karoui, Peng, Quenez~\cite{ElPeQu}).

\bl (Comparison Theorem) Given two coefficients $g_1$ and $g_2$
satisfying (A1) and (A2) and two terminal values $ \xi_1,\ \xi_2
\in L^{2}(\Omega, {\cal{F}}_{T}, P)$, we denote by $(y^1,z^1)$\
and $(y^2,z^2)$\ the solution of BSDE with the data $(\xi_1,g_1
)$\ and $(\xi_2,g_2 )$, respectively. Then we have:

{\rm (i) }(Monotonicity) If  $ \xi_1 \geq \xi_2$  and $ g_1 \geq
g_2, \ a.s.$, then $y^1_t\geq y^2_t,\ a.s.$, for all $t\in [0,
T].$

{\rm (ii)}(Strict Monotonicity) If, in addition to {\rm (i)}, we
also assume that $P(\xi_1 > \xi_2)> 0$, then $P\{y^1_t> y^2_t\}>0,
\ 0 \leq t \leq T,$\ and in particular, $ y^1_0> y^2_0.$ \el

Using the notation introduced in Lemma 2.2 we now suppose that,
for some $g: \Omega\times[0, T]\times{\mathbb{R}}
\times{\mathbb{R}}^{d}\longrightarrow {\mathbb{R}}$\ satisfying
(A1) and (A2) and for some $i\in \{1, 2\}$, the drivers $g_i, \
i=1, 2,$\ are of the form
$$g_i(s, y_s^i, z_s^i)=g(s, y_s^i, z_s^i)+\varphi_i(s),\ \ \mbox{dsdP-a.e.},\ i=1, 2,$$
where $\varphi_i\in {\cal{H}}^{2}(0,T;{\mathbb{R}}),\ i=1, 2.$\
Then, for terminal values $\xi_1,\ \xi_2\ \mbox{belonging to}\
L^{2}(\Omega, {\cal{F}}_{T}, P)$\ we have the following

 \bl The difference of the solutions $(y^1, z^1)$ and $(y^2, z^2)$ of BSDE with the data $(\xi_1, g_1)$\ and $(\xi_2, g_2)$, respectively, satisfies
 the following estimate:
 $$
  \begin{array}{ll}
  &|y^1_t-y^2_t|^2+\frac{1}{2}E[\int^T_te^{\beta(s-t)}[|
  y^1_s-y^2_s|^2+ |
  z^1_s-z^2_s|^2]ds|{\cal{F}}_t]
     \\
  \leq& E[e^{\beta(T-t)}|\xi_1-\xi_2|^2|{\cal{F}}_t]+ E[\int^T_te^{\beta(s-t)}
           |\varphi_1(s)-\varphi_2(s)|^2ds|{\cal{F}}_t],\ \mbox{P-a.s.,\ for all}\ 0\leq t\leq T,
  \end{array}
  $$
where $\beta=16(1+C^2)$.
  \el
For the proof the reader is referred to El Karoui, Peng,
Quenez~\cite{ElPeQu} or Peng~\cite{Pe1}.
\section{\large{Forward- Backward SDES (FBSDEs) }}

 \hskip1cm In this section we give an overview over basic results on BSDEs associated
 with Forward SDEs (for short: FSDEs). We consider measurable functions $b:[0,T]\times \Omega\times
{\mathbb{R}}^n\rightarrow {\mathbb{R}}^n \ $ and
         $\sigma:[0,T]\times \Omega\times {\mathbb{R}}^n\rightarrow {\mathbb{R}}^{n\times d}$
which are supposed to satisfy the following conditions:
 $$
  \begin{array}{ll}
\mbox{(i)}&b(\cdot,0)\ \mbox{and}\ \sigma(\cdot,0)\ \mbox{are} \
{\cal{F}}_t-\mbox{adapted processes, and there exists some}\\
 & \mbox{constant}\ C>0\  \mbox{such that}\\
 &\hskip 1cm|b(t,x)|+|\sigma(t,x)|\leq C(1+|x|), a.s.,\
                                  \mbox{for all}\ 0\leq t\leq T,\ x\in {\mathbb{R}}^n;\\
\mbox{(ii)}&b\ \mbox{and}\ \sigma\ \mbox{are Lipschitz in}\ x,\ \mbox{i.e., there is some constant}\ C>0\ \mbox{such that}\\
           &\hskip 1cm|b(t,x)-b(t,x')|+|\sigma(t,x)-\sigma(t,x')|\leq C| x-x'|,\ a.s.,\\
 & \hbox{ \ \ }\hskip7cm\mbox{for all}\ 0\leq t \leq T,\ x,\ x'\in {\mathbb{R}}^n.\\
 \end{array}
  \eqno{\mbox{(H3.1)}}
  $$\par
  We now consider the following SDE parameterized by the
  initial condition $(t,\zeta)\in[0,T]\times L^2(\Omega,{\cal{F}}_t,P;{\mathbb{R}}^n)$:
  \be
  \left\{
  \begin{array}{rcl}
  dX_s^{t,\zeta}&=&b(s,X_s^{t,\zeta})ds+\sigma(s,X_s^{t,\zeta})dB_s,\ s\in[t,T],\\
  X_t^{t,\zeta}&=&\zeta.
  \end{array}
  \right.
  \ee
Under the assumption (H3.1), SDE (3.1) has a unique strong solution
and, for any $p\geq 2,$\ there exists $C_{p}\in {\mathbb{R}}$\ such
that, for any $t\in[0,T]\ \mbox{and}\ \zeta,\zeta'\in
L^p(\Omega,{\cal{F}}_t,P;{\mathbb{R}}^n),$
 \be
 \begin{array}{rcl}
 E[\sup\limits_{t\leq s\leq T}| X_s^{t,\zeta}-X_s^{t,\zeta'}|^p|{\cal{F}}_t]
                             &\leq& C_{p}|\zeta-\zeta'|^p, \ \ a.s.,\\
  E[\sup\limits_{t\leq s\leq T}| X_s^{t,\zeta}|^p|{\cal{F}}_t]
                       &\leq& C_{p}(1+|\zeta|^p),\ \  a.s..
 \end{array}
\ee We emphasize that the constant $C_{p}$ in (3.2) only depends
on the Lipschitz and the growth constants of $b$ and $\sigma$.
 Let now be given two real valued functions $f(t,x,y,z)$ and $\Phi(x)$ which shall satisfy the
following conditions:
$$
\begin{array}{ll}
\mbox{(i)}&\Phi:\Omega\times {\mathbb{R}}^n\rightarrow
{\mathbb{R}} \ \mbox{is an}\
{\cal{F}}_T\otimes{\cal{B}}({\mathbb{R}}^n)
             \mbox{-measurable random variable and}\\
          & f:[0,T]\times \Omega\times {\mathbb{R}}^n\times {\mathbb{R}}\times
          {\mathbb{R}}^d \rightarrow {\mathbb{R}}\ \mbox{is a measurable process}\ \mbox{such that} \\
          &f(\cdot,x,y,z)\ \mbox{is}\ {\cal{F}}_t \mbox{-adapted, for all $(x, y, z)\in{\mathbb{R}}^n\times {\mathbb{R}}\times
          {\mathbb{R}}^d $;}\\
\mbox{(ii)}&\mbox{There exists a constant}\ C>0\ \mbox{such that}\\
          &| f(t,x,y,z)-f(t,x',y',z')| +| \Phi(x)-\Phi(x')|\\
          &\hskip 5cm\leq C(|x-x'|+ |y-y'|+|z-z'|),\ \  a.s.,\\
&\hskip 3cm\mbox{for all}\ 0\leq t\leq T,\ x,\ x'\in {\mathbb{R}}^n,\ y,\ y'\in {\mathbb{R}}\ \mbox{and}\ z,\ z'\in {\mathbb{R}}^d;\\
\mbox{(iii)}&f\ \mbox{and}\ \Phi \ \mbox{satisfy a linear growth condition, i.e., there exists some}\ C>0\\\
    & \mbox{such that, dt}\times \mbox{dP-a.e.},\ \mbox{for all}\ x\in
    {\mathbb{R}}^n,\\
    &\hskip 2cm|f(t,x,0,0)| + |\Phi(x)| \leq C(1+|x|).\\
\end{array}
\eqno{\mbox{(H3.2)}} $$
 With the help of the above assumptions we can verify that the coefficient $f(s,X_s^{t,\zeta},y,z)$\ satisfies the hypotheses (A1),
 (A2) and  $\xi=\Phi(X_T^{t,\zeta})$ $\in
 L^2(\Omega,{\cal{F}}_T,P;{\mathbb{R}})$. Therefore, the following BSDE
 possesses a unique solution:
\be
\left\{
\begin{array}{rcl}
-dY_s^{t,\zeta}&=&f(s,X_s^{t,\zeta},Y_s^{t,\zeta},Z_s^{t,\zeta})ds-Z_s^{t,\zeta}dB_s,\ s\in [t, T],\\
Y_T^{t,\zeta}&=&\Phi(X_T^{t,\zeta}).\\
\end{array}
\right.
\ee

\bp  We suppose that the hypotheses (H3.1) and (H3.2) hold. Then,
for any $0\leq t\leq T$ and the associated initial conditions
 $\zeta,\zeta'\in L^2(\Omega,{\cal{F}}_t,P;{\mathbb{R}}^n)$, we
have the following estimates:\\  $\mbox{}\hskip3cm\mbox{\rm(i)}
E[\sup\limits_{t\leq s\leq T}|Y_s^{t,\zeta}|^2 +
\int_t^T|Z_s^{t,\zeta}|^2ds|{{\cal{F}}_t}]\leq C(1+|\zeta|^2),\
a.s.; $\\
$\mbox{}\hskip3cm\mbox{\rm(ii)}E[\sup\limits_{t\leq s\leq
T}|Y_s^{t,\zeta}-Y_s^{t,\zeta'}|^2+\int_t^T|Z_s^{t,\zeta}-Z_s^{t,\zeta'}|^2ds|{{\cal{F}}_t}]\leq
C|\zeta-\zeta'|^2,\  a.s.. $\\ In particular, \be
 \begin{array}{lll}
\mbox{\rm(iii)}&|Y_t^{t,\zeta}|\leq C(1+|\zeta|),\  a.s.; \hskip3cm\\
\mbox{\rm(iv)}&|Y_t^{t,\zeta}-Y_t^{t,\zeta'}|\leq C|\zeta-\zeta'|,\  a.s.,\hskip3cm \\
\end{array}
\ee where the constant $C>0$\ depends only on the Lipschitz and
the growth constants of $b$,\ $\sigma$, $f$\ and $\Phi$. \ep
 \noindent The proof can be found in Peng~\cite{Pe1}.
\vskip 0.3cm
 Let us now introduce the random field:
\be u(t,x)=Y_s^{t,x}|_{s=t},\ (t, x)\in [0,
T]\times{\mathbb{R}}^n, \ee where $Y^{t,x}$ is the solution of
BSDE (3.3) with $x \in {\mathbb{R}}^n$\ at the place of $\zeta\in
L^2(\Omega,{\cal{F}}_t,P;{\mathbb{R}}^n).$\\

As a consequence of Proposition 3.1 we have that, for all $t \in
[0, T] $, P-a.s.,
 \be
\begin{array}{ll}
\mbox{(i)}&| u(t,x)-u(t,y)| \leq C|x-y|,\ \mbox{for all}\ x, y\in {\mathbb{R}}^n;\\
\mbox{(ii)}&| u(t,x)|\leq C(1+|x|),\ \mbox{for all}\ x\in {\mathbb{R}}^n.\\
\end{array}
\ee
 \br\ In the general situation $u$ is an adapted random function, that is, for any $x\in
 {\mathbb{R}}^n,\ u(\cdot,x)$ is an ${\cal{F}}_t-$adapted real valued process. Indeed, recall that
  $b, \sigma, f\ \mbox{and}\ \Phi$ all are ${\cal{F}}_t$-adapted random functions. On the other hand, it is well
  known that, under the additional assumption that the
  functions $$ b, \sigma, f\ \mbox{and}\ \Phi\ \mbox{are
  deterministic,}\eqno{\mbox{(H3.3)}}  $$
also $u$ is a deterministic function of $(t,x)$. \er

 The random field $u$\ and $Y^{t,\zeta},\ (t, \zeta)\in [0,
T]\times L^2(\Omega,{\cal{F}}_t,P;{\mathbb{R}}^n),$\ are related
by the following theorem.
 \bt Under the assumptions (H3.1) and (H3.2), for any $t\in [0, T]$\ and $\zeta\in
L^2(\Omega,{\cal{F}}_t,P;{\mathbb{R}}^n),$\ we have \be
u(t,\zeta)=Y_t^{t,\zeta},\ \mbox{ P-a.s.}. \ee \et The proof of
Theorem 3.1 can be found in Peng~\cite{Pe1}, we give it for the
reader's convenience. It makes use of the following definition.

 \noindent
\bde For any t $\in [0, T]$, a sequence $\{A_i\}_{i=1}^{N}\subset
{\cal{F}}_t\ (\mbox{with}\ 1\leq N\leq \infty)$ is called a
partition of $(\Omega, {\cal{F}}_t)$\ if\ \
$\cup_{i=1}^{N}A_i=\Omega$\ and $ A_i\cap A_j=\phi, \
\mbox{whenever}\ i\neq j.$ \ede
 \noindent \textbf{Proof} (of Theorem 3.1): We first consider
the case where $\zeta$ is a simple random variable of the form\
\be\zeta=\sum\limits^N\limits_{i=1}x_i\textbf{1}_{A_i},
                        \ee
where$\{A_i\}^N_{i=1}$\ is a finite partition of $(\Omega,{\cal{F}}_t)$\ and $x_i\in {\mathbb{R}}^n$,\ for $1\leq i\leq N.$\\
For each $i$, we put $(X_s^i,Y_s^i,Z_s^i)\equiv
                 (X_s^{t,x_i},Y_s^{t,x_i},Z_s^{t,x_i}).$ Then $X^i$ is the solution of the SDE
$$
X^i_s =x_i +\int^s_t b(r,X^i_r)dr +\int^s_t \sigma (r,X^i_r)dB_r,\
s\in [t,T],
$$\\
 and $(Y^i,Z^i)$ is the solution of the associated BSDE
$$
Y^i_s =\Phi(X^i_T) +\int^T_s f(r,X^i_r,Y^i_r,Z^i_r)dr
   -\int^T_s Z^i_r dB_r,\ s\in [t,T].
$$
The above two equations are multiplied by $\textbf{1}_{A_i}$\ and
summed up with respect to $i$. Thus, taking into account that
$\sum\limits_i \varphi (x_i)\textbf{1}_{A_i}=\varphi
(\sum\limits_i x_i \textbf{1}_{A_i})$, we get
$$
\begin{array}{rcl}
 \sum\limits_{i=1}\limits^{N} \textbf{1}_{A_i} X^i_s &=&\sum\limits _{i=1}\limits^{N}
  x_i \textbf{1}_{A_i}+ \int^s_t b(r,\sum\limits _{i=1}\limits^{N} \textbf{1}_{A_i} X^i_r )dr
  +\int^s_t \sigma (r,\sum\limits_{i=1}\limits^{N} \textbf{1}_{A_i} X^i_r)dB_r
\end{array}
$$and
$$
\begin{array}{rcl}
\sum\limits _{i=1}\limits^{N}\textbf{1}_{A_i} Y^i_s & = &
\Phi(\sum\limits _{i=1}\limits^{N} \textbf{1}_{A_i}
X^i_T)+\int^T_s f(r,\sum \limits_{i=1}^{N} \textbf{1}_{A_i} X^i_r,
  \sum\limits _{i=1}^{N} \textbf{1}_{A_i} Y^i_r,
 \sum\limits _{i=1}^{N} \textbf{1}_{A_i} Z^i_r)dr \\
  & & -\int^T_s \sum\limits _{i=1}^{N} \textbf{1}_{A_i} Z^i_r
  dB_r.
\end{array}
$$
Then the strong uniqueness property of the solution of the SDE and
the BSDE yields
$$
X^{t,\zeta}_s =\sum \limits_{i=1}^{N} X^i_s \textbf{1}_{A_i},\
(Y^{t,\zeta}_s,Z^{t,\zeta}_s) =(\sum \limits_{i=1}^{N}
\textbf{1}_{A_i} Y^i_s, \sum \limits_{i=1}^{N} \textbf{1}_{A_i}
Z^i_s),\ s\in [t, T].
$$
Finally, from $u(t,x_i)=Y^i_t,\ 1\leq i\leq N$, we deduce that
$$
Y^{t,\zeta}_t=\sum \limits_{i=1}^{N}
Y^i_t\textbf{1}_{A_i}=\sum\limits_{i=1}^{N}u(t,x_i)
\textbf{1}_{A_i} =u(t,\sum \limits_{i=1}^{N} x_i \textbf{1}_{A_i})
=u(t,\zeta).
$$
Therefore, for simple random variables, we have the desired
result.

Given a general $\zeta\in L^2 (\Omega ,{\mathcal{F}}_t
,P;{\mathbb{R}}^n)$ we can choose a sequence of simple random
variables $\{\zeta_i\}$ which
 converges to $\zeta$ in $L^2(\Omega ,{\mathcal{F}}_t
,P;{\mathbb{R}}^n)$. Consequently, from the estimates (3.4), (3.6)
and the first step of the proof, we have
$$
\begin{array}{lrcl}
&E|Y^{t,\zeta_i}_t-Y^{t,\zeta}_t|^2&\leq&CE|\zeta_i -\zeta|^2\rightarrow 0,\ i\rightarrow\infty,\\
\mbox{ }\hskip1cm&
E|u(t,\zeta_i)-u(t,\zeta)|^2 &\leq& CE|\zeta_i -\zeta|^2 \rightarrow 0,\ i\rightarrow\infty,\\
\hbox{and}\hskip1cm& Y^{t,\zeta_i}_t&=& u(t,\zeta_i),\ i\geq 1.
\end{array}
$$
Then the proof is complete.\endpf\vskip 0.3cm \br Under (H3.1),
(H3.2) and (H3.3) we know $u(t, x)$\ is
 $\frac{1}{2}-$H\"{o}lder continuous in $t$: There exists a constant C such that,
  for every $x\in {\mathbb{R}}^n,\ t,\ t'\in [0, T]$,
  $$
|u(t, x)-u(t', x)|\leq C(1+|x|)|t-t'|^{\frac{1}{2}}.
  $$
This inequality can be proved with the help of Theorem 3.1. Since,
on the other hand, a similar result but in a more general setting
will be proved later (see Theorem 4.2) we don't give the proof
here.

For the case of random coefficients $b,\ \sigma,\ f$ and $\Phi$\
we can state the following property.
 \er
 \br\mbox{ }Let
us suppose in addition to the assumptions (H3.1) and (H3.2) that
$\sigma(\omega,t, \cdot)$\ and\ $b(\omega,t, \cdot) \mbox{are
continuously differentiable with Lipschitz derivative such that,
}$\ for some constant C,
$$\begin{array}{lrcl}
&|D_x\sigma(\omega, t, x)|+ |D_xb(\omega,t,  x)|\leq C,\
\mbox{dtdP-a.e.},\
\mbox{for all}\ x\in {\mathbb{R}}^n;\\
&D_x\sigma(\omega, t, \cdot),\ D_xb(\omega,t,  \cdot) \ \mbox{are Lipschitz, uniformly in}\ (\omega, t).\\
\end{array}
$$Then the random field $u(\omega, t, x):\ \Omega\times [0, T]\times {\mathbb{R}}^n\rightarrow
{\mathbb{R}}$\ possesses a continuous version.
  \er
\noindent The proof uses a standard argument based on the
properties of the stochastic flow associated with (3.1).

\section{\large{Stochastic
Differential Games and Associated Dynamic Programming Principles
}}

\hskip1cm Now we want to study the stochastic differential game. The
set of admissible control processes ${\mathcal{U}}$ (resp.,
${\mathcal{V}}$) for the first (resp., second) player is the set of
all U (resp., V)-valued ${\mathcal{F}}_t$-progressively measurable
processes. The control state spaces U and V are supposed to be
compact metric spaces.

For given admissible controls $u(\cdot)\in {\mathcal{U}}$ and
$v(\cdot)\in {\mathcal{V}}$, the according orbit which regards $t$
as the initial time and $\zeta \in L^2 (\Omega ,{\mathcal{F}}_t,
P;{\mathbb{R}}^n)$ as the initial state is defined by the solution
of the following SDE:
  \be
  \left \{
  \begin{array}{llll}
  dX^{t,\zeta ;u, v}_s & = & b(s,X^{t,\zeta; u,v}_s, u_s, v_s) ds + \sigma(s,X^{t,\zeta; u,v}_s, u_s, v_s) dB_s,\ s\in
   [t,T], \\
   X^{t,\zeta ;u, v}_t  & = & \zeta,
   \end{array}
   \right.
  \ee
where the mappings
  $$
  \begin{array}{llll}
  &   b:[0,T]\times {\mathbb{R}}^n\times U\times V \rightarrow {\mathbb{R}}^n \
  \mbox{and}\ \   \sigma: [0,T]\times {\mathbb{R}}^n\times U\times V\rightarrow {\mathbb{R}}^{n\times d} \\
     \end{array}
  $$
  satisfy the following conditions:
  $$
  \begin{array}{ll}
 \rm{(i)}& \mbox{For every fixed}\ x\in {\mathbb{R}}^n,\ b(., x, ., .)\ \mbox{and}\ \sigma(., x, ., .)
    \ \mbox{are continuous in}\ (t,u,v);\\
 \rm{(ii)}&\mbox{There exists a }C>0\ \mbox{such that, for all}\ t\in [0,T],\ x, x'\in {\mathbb{R}}^n,\ u \in U,\ v \in V, \\
   &\hskip1cm |b(t,x,u,v)-b(t,x',u ,v)|+ |\sigma(t,x,u,v)-\sigma(t,x',u, v)|\leq C|x-x'|.\\
  \end{array}
  \eqno{\mbox{(H4.1)}}
  $$

From (H4.1) we can get the global linear growth conditions of b
and $\sigma$, i.e., the existence of some $C>0$\ such that, for
all $0 \leq t \leq T,\ u\in U,\ v \in V,\  x\in {\mathbb{R}}^n $,
  \be
  |b(t,x,u,v)| +|\sigma (t,x,u,v)| \leq C(1+|x|).
  \ee
Obviously, under the above assumptions, for any $u(\cdot)\in
{\mathcal{U}}$ and $v(\cdot)\in {\mathcal{V}}$, SDE (4.1) has a
unique strong solution. Moreover, for any $p\geq 2$, there exists
$C_p\in \mathbb{R}$\ such that, for any $t \in [0,T]$,
$u(\cdot)\in {\mathcal{U}}, v(\cdot)\in {\mathcal{V}}$\ and $
\zeta, \zeta'\in L^2 (\Omega ,{\mathcal{F}}_t,P;{\mathbb{R}}^n),$\
we also have the following estimates, P-a.s.:\\
 \be
\begin{array}{rcl}
E[\sup \limits_{s\in [t,T]}|X^{t,\zeta; u, v}_s -X^{t,\zeta';u,
v}_s|^p|{{\mathcal{F}}_t}]
& \leq & C_p|\zeta -\zeta'|^p, \\
E[ \sup \limits_{s\in [t,T]} |X^{t,\zeta
;u,v}_s|^p|{{\mathcal{F}}_t}] & \leq &
                        C_p(1+|\zeta|^p).
\end{array}
\ee The constant $C_p$ depends only on the Lipschitz
 and the linear growth constants of $b$\ and $\sigma$
with respect to $x$.

Let now be given two functions
$$
\Phi: {\mathbb{R}}^n \rightarrow {\mathbb{R}},\ f:[0,T]\times
{\mathbb{R}}^n \times {\mathbb{R}} \times {\mathbb{R}}^d \times U
\times V \rightarrow {\mathbb{R}}
$$
that satisfy the following conditions:
$$
\begin{array}{ll}
\rm{(i)}& \mbox{For every fixed}\ (x, y, z)\in {\mathbb{R}}^n
\times {\mathbb{R}} \times {\mathbb{R}}^d , f(., x, y, z,.,.)\
\mbox{is continuous in}\ (t,u,v)\ \mbox{and}\\
&\mbox{there exists a constant}\ C>0 \ \mbox{such that, for all}\
t\in [0,T],\ x, x'\in {\mathbb{R}}^n,\ y, y'\in
{\mathbb{R}},\ z, z'\\
&\in {\mathbb{R}}^d,\ u \in U \ \mbox{and}\ v \in V,\\
&\hskip3cm\begin{array}{l}
|f(t,x,y,z,u,v)-f(t,x',y',z',u,v)| \\
\hskip3cm \leq C(|x-x'|+|y-y'| +|z-z'|);
\end{array}\\
\rm{(ii)}&\mbox{There is a constant}\ C>0 \ \mbox{such that, for
all}\ x, x'\in {\mathbb{R}}^n,\\
 &\mbox{  }\hskip3cm |\Phi (x) -\Phi (x')|\leq C|x-x'|.
 \end{array}
 \eqno {\mbox{(H4.2)}}
 $$
From (H4.2) we see that $f$\ and $\Phi$\ also satisfy the global
linear growth condition in $x$, i.e., there exists some $C>0$\
such that, for all $0 \leq t \leq T,\  u\in U,\ v \in V,\ x\in
{\mathbb{R}}^n $,
  \be
   |f(t,x,0,0,u,v)|+|\Phi (x)| \leq C(1+|x|).
   \ee
For any $u(\cdot) \in {\mathcal{U}}, $\ $v(\cdot) \in
{\mathcal{V}}$\ and $\zeta \in L^2
(\Omega,{\mathcal{F}}_t,P;{\mathbb{R}}^n)$,  the mappings $\xi:=
\Phi(X^{t,\zeta; u, v}_T)$ and $g(s,y,z):= f(s,X^{t,\zeta; u,
v}_s,y,z,u_s,v_s)$ satisfy the conditions of Lemma 2.1 on the
interval $[t, T]$. Therefore, there exists a unique solution to
the following BSDE:
      \be
   \left \{\begin{array}{rcl}
   -dY^{t,\zeta; u, v}_s & = & f(s,X^{t,\zeta; u, v}_s, Y^{t,\zeta; u, v}_s, Z^{t,\zeta; u, v}_s,
                              u_s, v_s) ds -Z^{t,\zeta; u, v}_s dB_s,\\
      Y^{t,\zeta; u, v}_T  & = & \Phi (X^{t,\zeta; u, v}_T),
   \end{array}\right.
   \ee
where $X^{t,\zeta; u, v}$\ is introduced by equation (4.1).

 Moreover, in analogy to Proposition 3.1, we can see that
there exists some constant $C>0$\ such that, for all $0 \leq t
\leq T,\ \zeta, \zeta' \in L^2(\Omega ,
{\mathcal{F}}_t,P;{\mathbb{R}}^n),\ u(\cdot) \in {\mathcal{U}}\
\mbox{and}\ v(\cdot) \in {\mathcal{V}},$\ P-a.s.,
 \be
\begin{array}{ll}
 {\rm(i)} & |Y^{t,\zeta; u, v}_t -Y^{t,\zeta'; u, v}_t| \leq C|\zeta -\zeta'|; \\
 {\rm(ii)} & |Y^{t,\zeta; u, v}_t| \leq C (1+|\zeta|). \\
\end{array}
\ee

 We now introduce the following subspaces of admissible
controls:

\noindent\bde\ An admissible control process $u=\{u_r, r\in [t,
s]\}$ (resp., $v=\{v_r, r\in [t, s]\}$) for Player I (resp., II)
on $[t, s] (t<s\leq T)$\ is an ${\mathcal{F}}_r$-progressively
measurable process taking values in U (resp., V). The set of all
admissible controls for Player I (resp., II) on $[t, s]$ is
denoted by\ ${\mathcal{U}}_{t, s}$\ (resp., ${\mathcal{V}}_{t,
s}).$\ We identify two processes $u$\ and $\bar{u}$\ in\
${\mathcal{U}}_{t, s}$\ and write $u\equiv \bar{u}\ \mbox{on}\ [t,
s],$\ if $P\{u=\bar{u}\ \mbox{a.e. in}\ [t, s]\}=1.$\ Similarly we
interpret $v\equiv \bar{v}\ \mbox{on}\ [t, s]$\ in
${\mathcal{V}}_{t, s}$. \ede
 Finally, we have still to define the admissible strategies for the
game.

\bde A nonanticipative strategy for Player I on $[t, s] (t<s\leq
T)$ is a mapping $\alpha: {\mathcal{V}}_{t, s}\longrightarrow
{\mathcal{U}}_{t, s}$ such that, for any
${\mathcal{F}}_r$-stopping time $S: \Omega\rightarrow [t, s]$\ and
any $ v_1, v_2 \in {\mathcal{V}}_{t, s}$\ with $ v_1\equiv v_2\
\mbox {on}\ \textbf{[\![}t, S\textbf{]\!]},$ it holds
$\alpha(v_1)\equiv \alpha(v_2)\ \mbox {on}\ \textbf{[\![}t,
S\textbf{]\!]}$.\ Nonanticipative strategies for Player II on $[t,
s]$, $\beta: {\mathcal{U}}_{t, s}\longrightarrow {\mathcal{V}}_{t,
s}$,  are defined similarly. The set of all nonanticipative
strategies $\alpha: {\mathcal{V}}_{t,s}\longrightarrow
{\mathcal{U}}_{t,s}$ for Player I on $[t, s]$ is denoted by
${\cal{A}}_{t,s}$. The set of all nonanticipative strategies
$\beta: {\mathcal{U}}_{t,s}\longrightarrow {\mathcal{V}}_{t,s}$
for Player II on $[t, s]$ is denoted by ${\cal{B}}_{t,s}$.
\\
(\mbox{Recall that}\ $\textbf{[\![}t,
S\textbf{]\!]}=\{(r,\omega)\in [0, T]\times \Omega, t\leq r\leq
S(\omega)\})$.\ede

 Given the control
processes $u(\cdot)\in {\mathcal{U}}_{t,T}$\ and $ v(\cdot)\in
{\mathcal{V}}_{t,T} $\ we introduce the following associated cost
functional
 \be
J(t, x; u, v):= Y^{t, x; u, v}_t,\ (t, x)\in [0, T]\times
{\mathbb{R}}^n,\ee where the process $Y^{t, x; u, v}$ is defined
by BSDE (4.5).

 \noindent Similarly to the
proof of Theorem 3.1 we can get that, for any $t\in[0, T]$\ and
$\zeta \in L^2 (\Omega ,{\mathcal{F}}_t ,P; {\mathbb{R}}^n)$,
 \be J(t, \zeta; u, v) = Y^{t,\zeta; u, v}_t,\
 \mbox{P-a.s.}.
\ee
 Being particularly interested in the case of a deterministic $\zeta$, i.e., $\zeta=x\in {\mathbb{R}}^n$,
 we define the lower value function of our stochastic
differential game \be W(t,x):= \mbox{essinf}_{\beta \in
{\cal{B}}_{t,T}}\mbox{esssup}_{u \in {\mathcal{U}}_{t,T}}J(t,x;
u,\beta(u)) \ee
 and its upper value function
  \be U(t,x):= \mbox{esssup}_{\alpha \in
{\cal{A}}_{t,T}}\mbox{essinf}_{v \in {\mathcal{V}}_{t,T}}J(t,x;
\alpha(v),v). \ee

\br (1) For the convenience of the reader we recall that, given a
family of real-valued random variables $\eta_{\alpha},\ \alpha\in
I,$ a random variable $\eta$\ is said to be
$\mbox{essinf}_{\alpha \in I}\eta_{\alpha}$, if\\
{\rm i)}\ $\eta\leq\eta_\alpha,\ \mbox{P-a.s.,\ for any}\ \
\alpha \in I$;\\
{\rm ii)}\ if there is another random variable $\xi$ such that\ \
$\xi\leq\eta_\alpha,\ \mbox{P-a.s.,\ for any}\ \alpha \in
I,\ \mbox{then}\ \xi\leq\\
 \mbox{  }\hskip0.2cm \ \eta,\ \mbox{P-a.s..}$\\
The random variable $\mbox{esssup}_{\alpha \in I}\eta_{\alpha}$\
can be introduced now by the relation $$\mbox{esssup}_{\alpha \in
I}\eta_{\alpha}=-\mbox{essinf}_{\alpha \in I}(-\eta_{\alpha}).$$
Finally, recall that $\mbox{essinf}_{\alpha \in
I}\eta_{\alpha}=\mbox{inf}_{n\geq 1}\eta_{\alpha_n}$\ for some
denumerable family $(\alpha_n)\subset I;$\ $\mbox{esssup}_{\alpha
\in I}\eta_{\alpha}$\ has the same property.\\
(2) Obviously, under the assumptions (H4.1)-(H4.2), the lower
value function $W(t,x)$\ as well as the upper value function
$U(t,x)$ are well-defined and a priori they both are bounded
${\mathcal{F}}_{t}$-measurable random variables. But it turns out
that $W(t,x)$\ and $U(t,x)$\ are even deterministic. Indeed,
concentrating on the study of the properties of $W(t, x)$\ (the
function U(t,x) can be analyzed in a same manner) we can state the
following:\er
 \bp For any $[t, x]\in [0, T]\times {\mathbb{R}}^n$,
we have $W(t,x)=E[W(t,x)]$, P-a.s.. Identifying $W(t,x)$ with its
deterministic version $E[W(t,x)]$\ we can consider $W:[0, T]\times
{\mathbb{R}}^n\longrightarrow {\mathbb{R}}$ as a deterministic
function.\ep

\br Recall that the fact that the lower and upper value functions
defined by Fleming and Souganidis [8] are deterministic is an
immediate consequence of their definition. Indeed, for a game over
the time interval $[t,T]$ only control processes which are
independent of the past ${\cal F}_t$ are considered as admissible,
and since the admissible strategies are supposed to associate
admissible control processes of one player with those of the other
player, all the associated cost functionals are independent of
${\cal F}_t$ and hence deterministic.\er

 \noindent \textbf{Proof}: Let $H$ denote the Cameron-Martin space
 of all absolutely continuous elements $h\in \Omega$\ whose Radon-Nikodym
 derivative  $\dot{h}$\ belongs to $L^2([0, T],{\mathbb{R}}^d).$\\
For any $h \in H$, we define the mapping $\tau_h\omega:=\omega+h,\
\omega\in \Omega. $\ Obviously, $\tau_h: \Omega\rightarrow\Omega$\
is a bijection and its law is given by
$P\circ[\tau_h]^{-1}=\exp\{\int^T_0\dot{h}_sdB_s-\frac{1}{2}\int^T_0|\dot{h}_s|^2ds\}P.$\
Let $(t, x)\in [0, T]\times {\mathbb{R}}^n$\ be arbitrarily fixed,
and put $H_t=\{h\in H|h(\cdot)=h(\cdot\wedge t)\}.$\ We split now
the proof in the following steps:
 \vskip0.1cm
\noindent $1^{st}$ step: For any $u\in {\mathcal{U}}_{t,T}, \ v\in
{\mathcal{V}}_{t,T},\ h \in H_t,\ J(t, x; u,v)(\tau_h)= J(t, x;
u(\tau_h),v(\tau_h)),\ \mbox{P-a.s..}$ \vskip0.1cm
 Indeed, we apply the Girsanov transformation to SDE(4.1) (with
 $\zeta=x$) and compare the obtained equation with the SDE
 obtained from (4.1) by substituting the transformed control
 processes $u(\tau_h), v(\tau_h)$\ for $u$\ and $v$. Then, from the uniqueness of the solution of
(4.1) we get $X_s^{t,x; u,v}(\tau_h)=X_s^{t,x;
u(\tau_h),v(\tau_h)},$ $ \mbox{for any}\ s\in [t, T],\
\mbox{P-a.s..}$\ Furthermore, by a similar Girsanov transformation
argument we get from the uniqueness of the solution of BSDE (4.5),
$$Y_s^{t,x; u,v}(\tau_h)=Y_s^{t,x; u(\tau_h),v(\tau_h)},\ \mbox{for
any}\ s\in [t, T],\ \mbox{P-a.s.,}$$
$$Z_s^{t,x; u,v}(\tau_h)=Z_s^{t,x; u(\tau_h),v(\tau_h)},\  \mbox{dsdP-a.e. on}\ [t, T]\times\Omega.$$
That means $$J(t, x; u,v)(\tau_h)= J(t, x; u(\tau_h),v(\tau_h)),\
\mbox{P-a.s..}$$
 \vskip0.1cm
\noindent $2^{nd}$ step: For $\beta\in {\cal{B}}_{t,T}, \ h \in
H_t,$\ let $\beta^h(u):=\beta(u(\tau_{-h}))(\tau_h),\ u\in
{\mathcal{U}}_{t,T}.$\ Then $\beta^h\in {\cal{B}}_{t,T}.$
\vskip0.1cm Obviously, $\beta^h$\ maps ${\mathcal{U}}_{t,T}$\ into
${\mathcal{V}}_{t,T}$.\ Moreover, this mapping is nonanticipating.
Indeed, let $S: \Omega\rightarrow [t, T]$\ be an
${\mathcal{F}}_{r}$-stopping time and $ u_1, u_2 \in
{\mathcal{U}}_{t, T}$\ with $ u_1\equiv u_2\ \mbox {on}\
\textbf{[\![}t, S\textbf{]\!]}.$\ Then, obviously, $
u_1(\tau_{-h})\equiv u_2(\tau_{-h})\ \mbox {on}\ \textbf{[\![}t,
S(\tau_{-h})\textbf{]\!]}$ (notice that $ S(\tau_{-h})\ \mbox{is
still a}$ stopping time), and because $\beta\in {\cal{B}}_{t,T}$\ we
have $\beta(u_1(\tau_{-h}))\equiv \beta(u_2(\tau_{-h}))\ $ $ \mbox
{on}\ \textbf{[\![}t, S(\tau_{-h})\textbf{]\!]}$. Therefore,
$$\beta^h(u_1)=\beta(u_1(\tau_{-h}))(\tau_h)\equiv \beta(u_2(\tau_{-h}))(\tau_h)=\beta^h(u_2)\ \mbox
{on}\ \textbf{[\![}t, S\textbf{]\!]}.$$ \vskip0.1cm
\noindent$3^{rd}$ step: For all $h\in H_t$\ and $\beta\in
{\mathcal{B}}_{t, T}$\ we have:
$$\{\mbox{esssup}_{u \in {\mathcal{U}}_{t,T}}J(t,x;
u,\beta(u))\}(\tau_h)=\mbox{esssup}_{u \in
{\mathcal{U}}_{t,T}}\{J(t,x; u,\beta(u))(\tau_h)\},\
\mbox{P-a.s.}.
$$

Indeed, with the notation $I(t,x,\beta):=\mbox{esssup}_{u \in
{\mathcal{U}}_{t,T}}J(t,x; u,\beta(u)),\ \beta\in
{\mathcal{B}}_{t, T},$\ we have \ $I(t,x,\beta)\geq J(t,x;
u,\beta(u)),$\ and thus $I(t,x,\beta)(\tau_h)\geq J(t,x;
u,\beta(u))(\tau_h), \ \mbox{P-a.s.,\ for}$ $\mbox{ all}\ u\in
{\mathcal{U}}_{t,T}.$\ On the other hand, for any random variable
$\zeta$\ satisfying $\zeta\geq J(t,x; u,\beta(u))(\tau_h),$\ and
hence also $\zeta(\tau_{-h})\geq J(t,x; u,\beta(u)), \
\mbox{P-a.s.,\ for}\ \mbox{ all}\ u\in {\mathcal{U}}_{t,T},$\ we
have\ $\zeta(\tau_{-h})\geq I(t,x,\beta), \ $ $ \mbox{P-a.s.,}$\
i.e., $\zeta\geq I(t,x,\beta)(\tau_{h}), \ \mbox{P-a.s..}$\
Consequently,
$$I(t,x,\beta)(\tau_{h})=\mbox{esssup}_{u \in
{\mathcal{U}}_{t,T}}\{J(t,x; u,\beta(u))(\tau_h)\},\
\mbox{P-a.s.}$$
 \vskip0.1cm
 \noindent$4^{th}$ step: $W(t,x)$\ is invariant with respect
 to the Girsanov transformation $\tau_h$, i.e.,
  $$W(t,x)(\tau_{h})=W(t,x), \ \mbox{P-a.s., for any}\ h\in H. $$

Indeed, similarly to the third step we can show that for all $h\in
H_t$,
$$\{\mbox{essinf}_{\beta \in
{\mathcal{B}}_{t,T}}I(t,x;\beta)\}(\tau_h)=\mbox{essinf}_{\beta
\in {\mathcal{B}}_{t,T}}\{I(t,x; \beta)(\tau_h)\},\ \mbox{P-a.s..}
$$\ Then, from the first step to the third step we have, for any $h\in H_t,$
 $$
   \begin{array}{rcl}
   W(t,x)(\tau_{h}) & = & \mbox{essinf}_{\beta \in
{\mathcal{B}}_{t,T}}\mbox{esssup}_{u \in
{\mathcal{U}}_{t,T}}\{J(t,x; u,\beta(u))(\tau_h)\}\\
       & = &  \mbox{essinf}_{\beta \in
{\mathcal{B}}_{t,T}}\mbox{esssup}_{u \in
{\mathcal{U}}_{t,T}}J(t,x; u(\tau_h),\beta^h(u(\tau_h))\\
& = &  \mbox{essinf}_{\beta \in
{\mathcal{B}}_{t,T}}\mbox{esssup}_{u \in
{\mathcal{U}}_{t,T}}J(t,x; u,\beta^h(u))\\
& = &  \mbox{essinf}_{\beta \in
{\mathcal{B}}_{t,T}}\mbox{esssup}_{u \in
{\mathcal{U}}_{t,T}}J(t,x; u,\beta(u))\\
& = &W(t,x),\ \mbox{P-a.s.,}
   \end{array}
$$
where we have used
$\{u(\tau_h)|u(\cdot)\in{\mathcal{U}}_{t,T}\}={\mathcal{U}}_{t,T},\
\{\beta^h|\beta \in {\mathcal{B}}_{t,T} \}={\mathcal{B}}_{t,T}$\
in order to obtain the both latter equalities. Therefore,\ for any
$h\in H_t,\ W(t,x)$ $(\tau_{h})= W(t,x),\ \mbox{P-a.s.,}$\ and
since $W(t,x)$\ is ${\mathcal{F}}_{t}$-measurable, we have this
relation even for all $ h\in H.$

 The result of the $4^{th}$ step
combined with the following auxiliary Lemma 4.1 completes the
proof.\endpf

\bl Let $\zeta$\ be a random variable defined over our classical
Wiener space $(\Omega, {\mathcal{F}}_T, P)$, such that
 $\zeta(\tau_{h})=\zeta,\ \mbox{P-a.s., for any}\ h\in H.$\ Then
$\zeta=E\zeta,\ \mbox{P-a.s..}$\el

 \noindent \textbf{Proof}:\ Let $h\in H$\ and $ A\in
 {\mathcal{B}}({\mathbb{R}}).$\ Then,
$$\begin{array}{lll}
&E[\textbf{1}_{\{\zeta\in
A\}}\exp\{\int^T_0\dot{h}_sdB_s-\frac{1}{2}\int^T_0|\dot{h}_s|^2ds\}]\\
=&E[\textbf{1}_{\{\zeta(\tau_{-h})\in
A\}}\exp\{\int^T_0\dot{h}_sdB_s-\frac{1}{2}\int^T_0|\dot{h}_s|^2ds\}]\\
=&E[\textbf{1}_{\{\zeta\in A\}}],
\end{array}
$$
from where we deduce that
$$E[\textbf{1}_{\{\zeta\in
A\}}\exp\{\int^T_0\dot{h}_sdB_s\}]=E[\textbf{1}_{\{\zeta\in
A\}}]E[\exp\{\int^T_0\dot{h}_sdB_s\}],
$$
i.e., $\mbox{for any}\ \varphi\in L^2([0, T]; {\mathbb{R}}^d),$
\be E[\textbf{1}_{\{\zeta\in
A\}}\exp\{\int^T_0\varphi_sdB_s\}]=E[\textbf{1}_{\{\zeta\in
A\}}]E[\exp\{\int^T_0\varphi_sdB_s\}]. \ee
 Consequently, taking into consideration the arbitrariness of $A\in
 {\mathcal{B}}({\mathbb{R}})$\ and of $\varphi\in L^2([0, T];
 {\mathbb{R}}^d),$\ it follows the independence of $\zeta$\ of $B$\ and hence of ${\mathcal{F}}_T$, but
this is only possible for deterministic $\zeta$.\endpf

The first property of the lower value function $W(t,x)$ which we
present is an immediate consequence of (4.6) and (4.9).

\bl\mbox{  }There exists a constant $C>0$\ such that, for all $ 0
\leq t \leq T,\ x, x'\in {\mathbb{R}}^n$,\be
\begin{array}{llll}
&{\rm(i)} & |W(t,x)-W(t,x')| \leq C|x-x'|;  \\
&{\rm(ii)} & |W(t,x)| \leq C(1+|x|).
\end{array}
\ee \el \endpf \vskip0.3cm We now discuss (the generalized)
dynamic programming principle (DPP) for our stochastic
differential game (4.1), (4.5) and (4.9). For this end we have to
 define the family of (backward) semigroups associated with BSDE
 (4.5). This notion of stochastic backward semigroups was first
 introduced by Peng [14] which was applied to study the DPP for
 stochastic control problems. Our approach adapts Peng's ideas to the framework of stochastic differential games.

 Given the initial data $(t,x)$, a positive number $\delta\leq T-t$, admissible control
 processes $u(\cdot) \in {\mathcal{U}}_{t, t+\delta},\ v(\cdot) \in {\mathcal{V}}_{t, t+\delta}$\ and a real-valued
 random variable $\eta \in L^2 (\Omega,
{\mathcal{F}}_{t+\delta},P;{\mathbb{R}})$, we put \be G^{t, x; u,
v}_{s,t+\delta} [\eta]:= \tilde{Y}_s^{t,x; u, v},\ \hskip0.5cm
s\in[t, t+\delta], \ee where the couple $(\tilde{Y}_s^{t,x;u, v},
\tilde{Z}_s^{t,x;u, v})_{t\leq s \leq t+\delta}$ is the solution
of the following BSDE with the time horizon $t+\delta$:
$$
\left \{\begin{array}{rcl}
 -d\tilde{Y}_s^{t,x;u, v} \!\!\!& = &\!\!\! f(s,X^{t,x;u, v}_s ,\tilde{Y}_s^{t,x;u, v}, \tilde{Z}_s^{t,x;u, v},
                           u_{s}, v_{s})ds \\
  && -\tilde{Z}_s^{t,x; u, v} dB_s , \hskip 1cm s\in [t,t+\delta],\\
 \tilde{Y}_{t+\delta}^{t,x; u, v}\!\!\! & =& \!\!\!\eta ,
\end{array}\right.
$$
and $X^{t,x;u, v}$\ is the solution of SDE (4.1). Then, obviously,
for the solution $(Y^{t,x;u, v}, Z^{t,x;u, v})$\ of BSDE (4.5) we
have \be G^{t,x;u, v}_{t,T} [\Phi (X^{t,x; u, v}_T)] =G^{t,x;u,
v}_{t,t+\delta} [Y^{t,x;u, v}_{t+\delta}]. \ee Moreover,
$$
\begin{array}{rcl}
 J(t,x;u, v)& = &Y_t^{t,x;u, v}=G^{t,x;u, v}_{t,T} [\Phi (X^{t,x; u, v}_T)]
  =G^{t,x;u,v}_{t,t+\delta} [Y^{t,x;u, v}_{t+\delta}]\\
  &=&G^{t,x;u,v}_{t,t+\delta} [J(t+\delta,X^{t,x;u, v}_{t+\delta};u, v)].
\end{array}
$$

\br When $f$\ is independent of $(y, z)$\ it holds that
$$G^{t,x;u,v}_{s,t+\delta}[\eta]=E[\eta + \int_s^{t+\delta}
f(r,X^{t,x;u, v}_r,u_{r}, v_{r})dr|{\cal{F}}_s],\ \ s\in [t,
t+\delta].$$

\er

 \bt\mbox{}Under the
assumptions (H4.1) and (H4.2), the lower value function $W(t,x)$
obeys the following
 DPP : For any $0\leq t<t+\delta \leq T,\ x\in {\mathbb{R}}^n,$
 \be
W(t,x) =\mbox{essinf}_{\beta \in {\mathcal{B}}_{t,
t+\delta}}\mbox{esssup}_{u \in {\mathcal{U}}_{t,
t+\delta}}G^{t,x;u,\beta(u)}_{t,t+\delta} [W(t+\delta,
X^{t,x;u,\beta(u)}_{t+\delta})].
 \ee
  \et
\noindent \textbf{Proof}: To simplify notations we put
$$W_\delta(t,x) =\mbox{essinf}_{\beta \in {\mathcal{B}}_{t,
t+\delta}}\mbox{esssup}_{u \in {\mathcal{U}}_{t,
t+\delta}}G^{t,x;u,\beta(u)}_{t,t+\delta} [W(t+\delta,
X^{t,x;u,\beta(u)}_{t+\delta})].$$ The proof that $W_\delta(t,x)$\
coincides with $W(t,x)$\ will be split into a sequel of lemmata
which all are supposed to satisfy (H4.1) and (H4.2).

\bl $W_\delta(t,x)$\ is deterministic.\el The proof of this lemma
uses the same ideas as that of Proposition 4.1 so that it can be
omitted here.\endpf

\bl$W_\delta(t,x)\leq W(t,x).$\el

\noindent\textbf{ Proof}: Let $\beta\in {\mathcal{B}}_{t, T}$\ be
arbitrarily fixed. Then, given a $u_2(\cdot)\in
{\mathcal{U}}_{t+\delta, T},$\ we define as follows the
restriction $\beta_1$\ of $\beta$\ to ${\mathcal{U}}_{t+\delta,
T}:$
$$\beta_1(u_1):=\beta(u_1\oplus u_2 )|_{[t,
t+\delta]},\ \mbox{ }\ u_1(\cdot)\in {\mathcal{U}}_{t, t+\delta},
$$
where $u_1\oplus u_2:=u_1\textbf{1}_{[t,
t+\delta]}+u_2\textbf{1}_{(t+\delta, T]}$\ extends $u_1(\cdot)$\
to an element of ${\mathcal{U}}_{t, T}$. It is easy to check that
$\beta_1\in {\mathcal{B}}_{t, t+\delta}.$\ Moreover, from the
nonanticipativity property of $\beta$\ we deduce that $\beta_1$\
is independent of the special choice of $u_2(\cdot)\in
{\mathcal{U}}_{t+\delta, T}.$\ Consequently, from the definition
of $W_\delta(t,x),$
 \be W_\delta(t,x)\leq \mbox{esssup}_{u_1
\in {\mathcal{U}}_{t,
t+\delta}}G^{t,x;u_1,\beta_1(u_1)}_{t,t+\delta} [W(t+\delta,
X^{t,x;u_1,\beta_1(u_1)}_{t+\delta})],\ \mbox{P-a.s..} \ee We use
the notation $I_\delta(t, x, u, v):=G^{t,x;u,v}_{t,t+\delta}
[W(t+\delta, X^{t,x;u,v}_{t+\delta})]$\ and notice that there
exists a sequence $\{u_i^1,\ i\geq 1\}\subset {\mathcal{U}}_{t,
t+\delta}$\ such that
$$I_\delta(t, x, \beta_1):=\mbox{esssup}_{u_1 \in {\mathcal{U}}_{t,
t+\delta}}I_\delta(t, x, u_1, \beta_1(u_1))=\mbox{sup}_{i\geq
1}I_\delta(t, x, u_i^1, \beta_1(u_i^1)),\ \ \mbox{P-a.s.}.$$ For
any $\varepsilon>0,$\ we put $\widetilde{\Gamma}_i:=\{I_\delta(t,
x, \beta_1)\leq I_\delta(t, x, u_i^1,
\beta_1(u_i^1))+\varepsilon\}\in {\mathcal{F}}_{t},\ i\geq 1.$\
Then $\Gamma_1:=\widetilde{\Gamma}_1,\
\Gamma_i:=\widetilde{\Gamma}_i\backslash(\cup^{i-1}_{l=1}\widetilde{\Gamma}_l)\in
{\mathcal{F}}_{t},\ i\geq 2,$\ form an $(\Omega,
{\mathcal{F}}_{t})$-partition, and $u^\varepsilon_1:=\sum_{i\geq
1}\textbf{1}_{\Gamma_i}u_i^1$\ belongs obviously to
${\mathcal{U}}_{t, t+\delta}.$\ Moreover, from the
nonanticipativity of $\beta_1$\ we have
$\beta_1(u^\varepsilon_1)=\sum_{i\geq
1}\textbf{1}_{\Gamma_i}\beta_1(u_i^1),$\ and from the uniqueness
of the solution of the FBSDE, we deduce that $I_\delta(t, x,
u^\varepsilon_1, \beta_1(u^\varepsilon_1))=\sum_{i\geq
1}\textbf{1}_{\Gamma_i}I_\delta(t, x, u_i^1, \beta_1(u_i^1)),\
\mbox{P-a.s..}$\ Hence, \be
\begin{array}{llll}
W_\delta(t,x)\leq I_\delta(t, x, \beta_1)&\leq &\sum_{i\geq
1}\textbf{1}_{\Gamma_i}I_\delta(t, x, u_i^1, \beta_1(u_i^1))
+\varepsilon=I_\delta(t, x, u^\varepsilon_1,
\beta_1(u^\varepsilon_1))+\varepsilon\\
&=& G^{t,x;u^\varepsilon_1, \beta_1(u^\varepsilon_1)}_{t,t+\delta}
[W(t+\delta, X^{t,x;u^\varepsilon_1,
\beta_1(u^\varepsilon_1)}_{t+\delta})]+\varepsilon,\
\mbox{P-a.s..}
\end{array}
\ee
 On the other hand, using the fact that $\beta_1(\cdot):=\beta(\cdot\oplus u_2
)\in {\mathcal{B}}_{t, t+\delta}$\ does not depend on $u_2(\cdot)\in
{\mathcal{U}}_{t+\delta, T}$\ we can define
$\beta_2(u_2):=\beta(u^\varepsilon_1\oplus u_2)|_{[t+\delta, T]},\
\mbox{for all }\ u_2(\cdot)\in {\mathcal{U}}_{t+\delta, T}. $\ The
such defined $\beta_2: {\mathcal{U}}_{t+\delta, T}\rightarrow
{\mathcal{V}}_{t+\delta, T}$\ belongs to ${\mathcal{B}}_{t+\delta,
T}\ \mbox{since}\ \beta\in {\mathcal{B}}_{t, T}$. Therefore, from
the definition of $W(t+\delta,y)$\ we have, for any $y\in
{\mathbb{R}}^n,$
$$W(t+\delta,y)\leq \mbox{esssup}_{u_2 \in {\mathcal{U}}_{t+\delta, T}}J(t+\delta, y; u_2, \beta_2(u_2)),\ \mbox{P-a.s..}$$
Finally, because there exists a constant $C\in {\mathbb{R}}$\ such
that \be
\begin{array}{llll}
{\rm(i)} & |W(t+\delta,y)-W(t+\delta,y')| \leq C|y-y'|,\ \mbox{for any}\ y,\ y' \in {\mathbb{R}}^n;  \\
{\rm(ii)} & |J(t+\delta, y, u_2, \beta_2(u_2))-J(t+\delta, y',
u_2, \beta_2(u_2))| \leq C|y-y'|,\ \mbox{P-a.s.,}\\
 &\mbox{ }\hskip1cm \mbox{for any}\ u_2\in {\mathcal{U}}_{t+\delta, T},
\end{array}
\ee (see Lemma 4.2-(i) and (4.6)-(i)) we can show by approximating
$X^{t,x;u_1^\varepsilon,\beta_1(u_1^\varepsilon)}_{t+\delta}$\
that
$$W(t+\delta, X^{t,x;u_1^\varepsilon,\beta_1(u_1^\varepsilon)}_{t+\delta} )\leq
\mbox{esssup}_{u_2 \in {\mathcal{U}}_{t+\delta, T}}J(t+\delta,
X^{t,x;u_1^\varepsilon,\beta_1(u_1^\varepsilon)}_{t+\delta}; u_2,
\beta_2(u_2)),\ \mbox{P-a.s..}$$ To estimate the right side of the
latter inequality we note that there exists some sequence $\{u_j^2,\
j\geq 1\}\subset {\mathcal{U}}_{t+\delta, T}$\ such that
$$\mbox{esssup}_{u_2 \in {\mathcal{U}}_{t+\delta,
T}}J(t+\delta,X^{t,x;u_1^\varepsilon,\beta_1(u_1^\varepsilon)}_{t+\delta};
u_2, \beta_2(u_2))=\mbox{sup}_{j\geq
1}J(t+\delta,X^{t,x;u_1^\varepsilon,\beta_1(u_1^\varepsilon)}_{t+\delta};
u^2_j, \beta_2(u^2_j)),\ \mbox{P-a.s..}$$
 Then, putting\\
$\widetilde{\Delta}_j:=\{\mbox{esssup}_{u_2 \in
{\mathcal{U}}_{t+\delta,
T}}J(t+\delta,X^{t,x;u_1^\varepsilon,\beta_1(u_1^\varepsilon)}_{t+\delta};
u_2, \beta_2(u_2))\leq
J(t+\delta,X^{t,x;u_1^\varepsilon,\beta_1(u_1^\varepsilon)}_{t+\delta};
u^2_j, \beta_2(u^2_j))+\varepsilon\}\in {\mathcal{F}}_{t+\delta},\
j\geq 1;$\ we have with $\Delta_1:=\widetilde{\Delta}_1,\
\Delta_j:=\widetilde{\Delta}_j\backslash(\cup^{j-1}_{l=1}\widetilde{\Delta}_l)\in
{\mathcal{F}}_{t+\delta},\ j\geq 2,$\ an $(\Omega,
{\mathcal{F}}_{t+\delta})$-partition and
$u^\varepsilon_2:=\sum_{j\geq 1}\textbf{1}_{\Delta_j}u_j^2$\
 $\in {\mathcal{U}}_{t+\delta, T}.$ From
the nonanticipativity of $\beta_2$\ we have
$\beta_2(u^\varepsilon_2)=\sum_{j\geq
1}\textbf{1}_{\Delta_j}\beta_2(u_j^2)$\ and from the definition of
$\beta_1,\ \beta_2$\ we know that $\beta(u_1^\varepsilon\oplus
u_2^\varepsilon)=\beta_1(u_1^\varepsilon)\oplus
\beta_2(u_2^\varepsilon ).$\ Thus, again from the uniqueness of
the solution of our FBSDE, we get
$$\begin{array}{lcl}
J(t+\delta,X^{t,x;u_1^\varepsilon,\beta_1(u_1^\varepsilon)}_{t+\delta};
u_2^\varepsilon,
\beta_2(u_2^\varepsilon))&=&Y_{t+\delta}^{t+\delta,X^{t,x;u_1^\varepsilon,\beta_1(u_1^\varepsilon)}_{t+\delta};
u_2^\varepsilon, \beta_2(u_2^\varepsilon)}\ \hskip2cm \mbox{(see (4.8))}\\
&=&\sum_{j\geq
1}\textbf{1}_{\Delta_j}Y_{t+\delta}^{t+\delta,X^{t,x;u_1^\varepsilon,\beta_1(u_1^\varepsilon)}_{t+\delta};
 u_j^2, \beta_2( u_j^2)}\\
&=&\sum_{j\geq
1}\textbf{1}_{\Delta_j}J(t+\delta,X^{t,x;u_1^\varepsilon,\beta_1(u_1^\varepsilon)}_{t+\delta};
u_j^2, \beta_2(u_j^2)),\ \mbox{P-a.s..}
\end{array}$$
Consequently, \be
\begin{array}{lll}
W(t+\delta,
X^{t,x;u_1^\varepsilon,\beta_1(u_1^\varepsilon)}_{t+\delta} )&\leq
& \mbox{esssup}_{u_2 \in {\mathcal{U}}_{t+\delta,
T}}J(t+\delta,X^{t,x;u_1^\varepsilon,\beta_1(u_1^\varepsilon)}_{t+\delta};
u_2, \beta_2(u_2))\\
&\leq& \sum_{j\geq
1}\textbf{1}_{\Delta_j}Y_{t+\delta}^{t,x;u_1^\varepsilon\oplus
u_j^2,\beta(u_1^\varepsilon\oplus
u_j^2)}+\varepsilon\\
& = & Y_{t+\delta}^{t,x;u_1^\varepsilon\oplus u^\varepsilon_2,
\beta(u_1^\varepsilon\oplus
u^\varepsilon_2)}+\varepsilon\\
&=&Y_{t+\delta}^{t,x;u^\varepsilon,\beta(u^\varepsilon)}+\varepsilon,\
 \mbox{P-a.s.,}
\end{array}
\ee where $u^\varepsilon:= u_1^\varepsilon\oplus
u^\varepsilon_2\in {\mathcal{U}}_{t, T}.$\ From (4.17), (4.19),
Lemma 2.2 (comparison theorem for BSDEs) and Lemma 2.3 we have
 \be
\begin{array}{lll}
W_\delta(t,x)&\leq& G^{t,x;u^\varepsilon_1,
\beta_1(u^\varepsilon_1)}_{t,t+\delta}
[Y_{t+\delta}^{t,x;u^\varepsilon,\beta(u^\varepsilon)}+\varepsilon]+\varepsilon \\
&\leq& G^{t,x;u^\varepsilon_1,
\beta_1(u^\varepsilon_1)}_{t,t+\delta}
[Y_{t+\delta}^{t,x;u^\varepsilon,\beta(u^\varepsilon)}]+
(C+1)\varepsilon\\
& =& G^{t,x;u^\varepsilon, \beta(u^\varepsilon)}_{t,t+\delta}
[Y_{t+\delta}^{t,x;u^\varepsilon,\beta(u^\varepsilon)}]+
(C+1)\varepsilon\\
& =& Y_{t}^{t,x;u^\varepsilon,\beta(u^\varepsilon)}+
(C+1)\varepsilon\\
&\leq& \mbox{esssup}_{u \in {\mathcal{U}}_{t,
T}}Y_{t}^{t,x;u,\beta(u)}+ (C+1)\varepsilon,\ \mbox{P-a.s..}
\end{array}
\ee Since $\beta\in {\mathcal{B}}_{t, T}$\ has been arbitrarily
chosen we have (4.20) for all $\beta\in {\mathcal{B}}_{t, T}$.
Therefore,
 \be W_\delta(t,x)\leq \mbox{essinf}_{\beta\in
{\mathcal{B}}_{t, T}}\mbox{esssup}_{u \in {\mathcal{U}}_{t,
T}}Y_{t}^{t,x;u,\beta(u)}+ (C+1)\varepsilon= W(t, x)+
(C+1)\varepsilon.\ee Finally, letting $\varepsilon\downarrow0,\
\mbox{we get}\ W_\delta(t,x)\leq W(t, x).$\endpf

\bl$ W(t, x)\leq W_\delta(t,x).$\el

 \noindent \textbf{Proof}: We continue to use the notations introduced above, from the definition of
$W_\delta(t,x)$\ we have
$$
\begin{array}{lll}
W_\delta(t,x)&=& \mbox{essinf}_{\beta_1 \in {\mathcal{B}}_{t,
t+\delta}}\mbox{esssup}_{u_1 \in {\mathcal{U}}_{t,
t+\delta}}G^{t,x;u_1,\beta_1(u_1)}_{t,t+\delta} [W(t+\delta,
X^{t,x;u_1,\beta_1(u_1)}_{t+\delta})]\\
&=&\mbox{essinf}_{\beta_1 \in {\mathcal{B}}_{t,
t+\delta}}I_\delta(t, x, \beta_1),
\end{array}
$$
and, for some sequence $\{\beta_i^1,\ i\geq 1\}\subset
{\mathcal{B}}_{t, t+\delta},$
$$W_\delta(t,x)=\mbox{inf}_{i\geq
1}I_\delta(t, x, \beta_i^1),\ \mbox{P-a.s..}$$ For any
$\varepsilon>0,$\ we let $\widetilde{\Lambda}_i:=\{I_\delta(t, x,
\beta_i^1)-\varepsilon\leq W_\delta(t,x)\}\in {\mathcal{F}}_{t},\
i\geq 1,$ $\Lambda_1:=\widetilde{\Lambda}_1\ \mbox{and}\
\Lambda_i:=\widetilde{\Lambda}_i\backslash(\cup^{i-1}_{l=1}\widetilde{\Lambda}_l)\in
{\mathcal{F}}_{t},\ i\geq 2.$\ Then $\{\Lambda_i,\ i\geq 1\}$\ is
an $(\Omega, {\mathcal{F}}_{t})$-partition,
$\beta^\varepsilon_1:=\sum_{i\geq
1}\textbf{1}_{\Lambda_i}\beta_i^1$\ belongs to ${\mathcal{B}}_{t,
t+\delta},$\ and from the uniqueness of the solution of our FBSDE
we conclude that $I_\delta(t, x, u_1,
\beta^\varepsilon_1(u_1))=\sum_{i\geq
1}\textbf{1}_{\Lambda_i}I_\delta(t, x, u_1, \beta_i^1(u_1)),\
\mbox{P-a.s., for all}$\ \ $u_1(\cdot)\in {\mathcal{U}}_{t,
t+\delta}.$\ Hence,
 \be
\begin{array}{lll}
W_\delta(t,x)&\geq &\sum_{i\geq
1}\textbf{1}_{\Lambda_i}I_\delta(t, x,
\beta_i^1)-\varepsilon\\
&\geq&\sum_{i\geq 1}\textbf{1}_{\Lambda_i}I_\delta(t, x, u_1,
\beta_i^1(u_1))-\varepsilon \\
&=& I_\delta(t, x, u_1, \beta^\varepsilon_1(u_1))-\varepsilon\\
&=& G^{t,x;u_1, \beta^\varepsilon_1(u_1)}_{t,t+\delta}
[W(t+\delta, X^{t,x;u_1,
\beta_1^\varepsilon(u_1)}_{t+\delta})]-\varepsilon,\ \mbox{P-a.s.,
for all}\ \ u_1\in {\mathcal{U}}_{t, t+\delta}.
\end{array}
\ee
 On the other hand, from the definition of $W(t+\delta,y),$\
with the same technique as before, we deduce that, for any $y\in
{\mathbb{R}}^n,$\ there exists $\beta^\varepsilon_y\in
{\mathcal{B}}_{t+\delta, T}$\ \ such that \be W(t+\delta,y)\geq
\mbox{esssup}_{u_2 \in {\mathcal{U}}_{t+\delta,T}}J(t+\delta, y;
u_2, \beta^\varepsilon_y(u_2))-\varepsilon,\ \mbox{P-a.s..}\ee Let
$\{O_i\}_{i\geq1}\subset {\mathcal{B}}({\mathbb{R}}^n)$\ be a
decomposition of ${\mathbb{R}}^n$\ such that
$\sum\limits_{i\geq1}O_i={\mathbb{R}}^n\ \mbox{and}\
\mbox{diam}(O_i)\leq \varepsilon,\ i\geq 1.$\ And let $y_i$\ be an
arbitrarily fixed element of $O_i,\ i\geq1.$\ Defining
$[X^{t,x;u_1,
\beta_1^\varepsilon(u_1)}_{t+\delta}]:=\sum\limits_{i\geq1}y_i\textbf{1}_{\{X^{t,x;u_1,
\beta_1^\varepsilon(u_1)}_{t+\delta}\in O_i\}},$\ we have \be
|X^{t,x;u_1, \beta_1^\varepsilon(u_1)}_{t+\delta}-[X^{t,x;u_1,
\beta_1^\varepsilon(u_1)}_{t+\delta}]|\leq \varepsilon,\
\mbox{everywhere on}\ \Omega, \ \mbox{for all}\ u_1\in
{\mathcal{U}}_{t, t+\delta}.\ee\ Moreover, for each $y_i,$\ there
exists some $\beta^\varepsilon_{y_i}\in {\mathcal{B}}_{t+\delta,
T}$\ such that (4.23) holds, and, clearly,
$\beta^{\varepsilon}_{u_1}:=\sum\limits_{i\geq1}\textbf{1}_{\{X^{t,x;u_1,
\beta_1^\varepsilon(u_1)}_{t+\delta}\in
O_i\}}\beta^\varepsilon_{y_i}\in {\mathcal{B}}_{t+\delta, T}.$

 Now we can define the new strategy
$\beta^{\varepsilon}(u):=\beta_1^\varepsilon(u_1)\oplus
\beta^{\varepsilon}_{u_1}(u_2),\ u\in {\mathcal{U}}_{t, T},\
\mbox{where}\ u_1=u|_{[t, t+\delta]},\ u_2=u|_{(t+\delta, T]}$\
(restriction of $u$ to $[t, t+\delta]\times \Omega$\ and
$(t+\delta, T]\times \Omega$, resp.). Obviously,
$\beta^{\varepsilon}$\ maps ${\mathcal{U}}_{t,T}$\ into
${\mathcal{V}}_{t,T}.$\ Moreover, $\beta^{\varepsilon}$\ is
nonanticipating: Indeed, let $S: \Omega\longrightarrow[t, T]$\ be
an ${\mathcal{F}}_r$-stopping time and $u, u'\in
{\mathcal{U}}_{t,T}$\ be such that $u\equiv u'$\ on
$\textbf{[\![}t, S\textbf{]\!]}$. Decomposing $u,\ u'$\ into $u_1,
u'_1\in {\mathcal{U}}_{t,t+\delta},\ u_2, u'_2\in
{\mathcal{U}}_{t+\delta, T}$\ such that $u=u_1\oplus u_2\
\mbox{and}\ u'=u'_1\oplus u'_2.$\ We have $u_1\equiv u_1'$\ on
$\textbf{[\![}t, S\wedge(t+\delta)\textbf{]\!]}$\ from where we
get $\beta_1^\varepsilon(u_1)\equiv \beta_1^\varepsilon(u_1')$\ on
$\textbf{[\![}t, S\wedge(t+\delta)\textbf{]\!]}$\ (recall that
$\beta_1^\varepsilon$\ is nonanticipating). On the other hand,
$u_2\equiv u_2'$\ on $\textbf{]\!]}t+\delta,
S\vee(t+\delta)\textbf{]\!]}(\subset (t+\delta,T]\times
\{S>t+\delta\}),$\ and on $\{S>t+\delta\}$\ we have $X^{t,x;u_1,
\beta_1^\varepsilon(u_1)}_{t+\delta}=X^{t,x;u'_1,
\beta_1^\varepsilon(u'_1)}_{t+\delta}.$\ Consequently, from our
definition,
$\beta^{\varepsilon}_{u_1}=\beta^{\varepsilon}_{u'_1}$\ on
$\{S>t+\delta\}$\ and
$\beta^{\varepsilon}_{u_1}(u_2)\equiv\beta^{\varepsilon}_{u'_1}(u'_2)$\
on $\textbf{]\!]}t+\delta, S\vee(t+\delta)\textbf{]\!]}.$ This
yields $\beta^{\varepsilon}(u)=\beta_1^\varepsilon(u_1)\oplus
\beta^{\varepsilon}_{u_1}(u_2)\equiv\beta_1^\varepsilon(u'_1)\oplus
\beta^{\varepsilon}_{u'_1}(u'_2)=\beta^{\varepsilon}(u')$\ on
$\textbf{[\![}t, S\textbf{]\!]}$, from where it follows that
$\beta^{\varepsilon}\in {\mathcal{B}}_{t, T}.$

Let now $u\in {\mathcal{U}}_{t, T}$\ be arbitrarily chosen and
decomposed into $u_1=u|_{[t, t+\delta]}\in {\mathcal{U}}_{t,
t+\delta}$\ and $u_2=u|_{(t+\delta, T]}\in
{\mathcal{U}}_{t+\delta, T}.$\ Then, from (4.22), (4.18)-(i),
(4.24) and the lemmata 2.2 (comparison theorem) and 2.3 we obtain,
\be
\begin{array}{llll}
W_\delta(t,x)&\geq&  G^{t,x;u_1,
\beta^\varepsilon_1(u_1)}_{t,t+\delta} [W(t+\delta, X^{t,x;u_1,
\beta_1^\varepsilon(u_1)}_{t+\delta})]-\varepsilon\\
&\geq&G^{t,x;u_1,
\beta^\varepsilon_1(u_1)}_{t,t+\delta}[W(t+\delta, [X^{t,x;u_1,
\beta_1^\varepsilon(u_1)}_{t+\delta}])-C\varepsilon]-\varepsilon\\
&\geq&G^{t,x;u_1,
\beta^\varepsilon_1(u_1)}_{t,t+\delta}[W(t+\delta, [X^{t,x;u_1,
\beta_1^\varepsilon(u_1)}_{t+\delta}])]- C\varepsilon\\
&=&G^{t,x;u_1,\beta_1^\varepsilon(u_1)}_{t,t+\delta}[\sum\limits_{i\geq1}\textbf{1}_{\{X^{t,x;u_1,
\beta_1^\varepsilon(u_1)}_{t+\delta}\in O_i\}}W(t+\delta,y_i)]-
C\varepsilon,\ \ \mbox{P-a.s.}.
\end{array}
\ee Furthermore, from (4.23), (4.18)-(ii), (4.24), Lemmata 2.2
(comparison theorem) and 2.3, we have, \be
\begin{array}{lcl}
W_\delta(t,x)&\geq&
G^{t,x;u_1,\beta_1^\varepsilon(u_1)}_{t,t+\delta}[\sum\limits_{i\geq1}\textbf{1}_{\{X^{t,x;u_1,
\beta_1^\varepsilon(u_1)}_{t+\delta}\in O_i\}}J(t+\delta, y_i;
u_2, \beta^\varepsilon_{y_i}(u_2))-\varepsilon]-
C\varepsilon\\
&\geq&G^{t,x;u_1,\beta_1^\varepsilon(u_1)}_{t,t+\delta}[\sum\limits_{i\geq1}\textbf{1}_{\{X^{t,x;u_1,
\beta_1^\varepsilon(u_1)}_{t+\delta}\in O_i\}}J(t+\delta, y_i;
u_2, \beta^\varepsilon_{y_i}(u_2))]-C\varepsilon\\
&=&G^{t,x;u_1,\beta_1^\varepsilon(u_1)}_{t,t+\delta}[J(t+\delta,
[X^{t,x;u_1, \beta_1^\varepsilon(u_1)}_{t+\delta}]; u_2,
\beta^\varepsilon_{u_1}(u_2))]-C\varepsilon\\
&\geq &G^{t,x;u_1,
\beta^\varepsilon_1(u_1)}_{t,t+\delta}[J(t+\delta,X^{t,x;u_1,
\beta_1^\varepsilon(u_1)}_{t+\delta}; u_2,
\beta^\varepsilon_{u_1}(u_2))-C\varepsilon]- C\varepsilon\\
&\geq &G^{t,x;u_1,
\beta^\varepsilon_1(u_1)}_{t,t+\delta}[J(t+\delta,X^{t,x;u_1,
\beta_1^\varepsilon(u_1)}_{t+\delta}; u_2,
\beta^\varepsilon_{u_1}(u_2))]- C\varepsilon\\
&=& G^{t,x;u,\beta^\varepsilon(u)}_{t,t+\delta}[Y_{t+\delta}^{t,
x, u, \beta^{\varepsilon}(u)}]- C\varepsilon\\
&=& Y_{t}^{t, x; u, \beta^{\varepsilon}(u)}- C\varepsilon,\
\mbox{P-a.s., for any}\ u\in {\mathcal{U}}_{t, T}.
\end{array}
\ee Consequently, \be
\begin{array}{llll}
W_\delta(t,x)&\geq& \mbox{esssup}_{u \in {\mathcal{U}}_{t,
T}}J(t, x; u, \beta^{\varepsilon}(u))- C\varepsilon\\
&\geq&\mbox{essinf}_{\beta \in {\mathcal{B}}_{t,
T}}\mbox{esssup}_{u \in {\mathcal{U}}_{t, T}}J(t, x; u,
\beta(u))- C\varepsilon\\
&=&W(t,x)- C\varepsilon,\ \mbox{P-a.s.}.
\end{array}
\ee Finally, letting $\varepsilon\downarrow0$\ we get
$W_\delta(t,x)\geq W(t,x).$\ The proof is complete.\endpf

\br\mbox{}{\rm{(i)}} From the inequalities (4.17) and (4.22) we
see that for all $(t, x)\in [0,T]\times {\mathbb{R}}^n,$\
$\delta>0$\ with $0<\delta\leq T-t$\ and $\varepsilon>0$,\ it
holds:\\
 a) For every $\beta \in {\cal{B}}_{t, t+\delta},$\
there exists some $u^{\varepsilon}(\cdot) \in {\cal{U}}_{t,
t+\delta}$\ such that
 \be W(t,x)(=W_\delta(t, x))\leq G^{t,x;
u^{\varepsilon},\beta(u^{\varepsilon})}_{t,t+\delta}
      [W(t+\delta, X^{t,x; u^{\varepsilon},\beta(u^{\varepsilon})}_{t+\delta})]+
      \varepsilon,\ \mbox{P-a.s..}
\ee
 b) There exists some $\beta^{\varepsilon} \in {\cal{B}}_{t,
t+\delta}$\ such that, for all $u\in {\cal{U}}_{t, t+\delta},$ \be
W(t,x)(=W_\delta(t, x))\geq G^{t,x;
u,\beta^{\varepsilon}(u)}_{t,t+\delta}
      [W(t+\delta, X^{t,x;u,\beta^{\varepsilon}(u)}_{t+\delta})]-
      \varepsilon,\ \mbox{P-a.s.}.
\ee {\rm{(ii)}} Recall that the lower value function $W$\ is
deterministic. Thus, by taking the expectation on both sides of
(4.28) and (4.29) we can show that $$ W(t,x)= \mbox{inf}_{\beta
\in {\cal{B}}_{t,T}}\mbox{sup}_{u \in {\mathcal{U}}_{t,T}}E[J(t,x;
u,\beta(u))]. $$ In analogy we also have
$$ U(t,x)= \mbox{sup}_{\alpha
\in {\cal{A}}_{t,T}}\mbox{inf}_{v \in {\mathcal{V}}_{t,T}}E[J(t,x;
\alpha(v), v)]. $$ The above formulas look similar to the
definitions of the lower and the upper value functions defined by
Fleming and Souganidis~\cite{FS1} for the case of $f$\ being
independent of $(y, z)$. However, they consider only control
processes which are independent of the past ${\mathcal{F}}_{t}$. In
Remark 6.3 we will come back to this comparison and identify their
value functions with ours for such coefficient $f$.\er

In Lemma 4.2 we have already seen that the lower value function
$W(t,x)$\ is Lipschitz continuous in $x$, uniformly in $t$. With
the help of Theorem 4.1 we can now also study the continuity
properties of $W(t,x)$\ in $t$.
 \bt\mbox{ }Let us suppose that the assumptions (H4.1) and (H4.2)
hold. Then the lower value function $W(t,x)$ is
 $\frac{1}{2}-$H\"{o}lder continuous in $t$: There exists a constant C such that,
  for every $x\in {\mathbb{R}}^n,\ t, t'\in [0, T]$,
  $$
|W(t, x)-W(t', x)|\leq C(1+|x|)|t-t'|^{\frac{1}{2}}.
  $$
  \et
\noindent \textbf{Proof}:  Let $(t, x)\in [0,T]\times
{\mathbb{R}}^n$\ and $\delta>0$\ be arbitrarily given such that
$0<\delta\leq T-t$. Our objective is to prove the following
inequality by using (4.28) and (4.29):
 \be -C(1+|x|)\delta^{\frac{1}{2}}\leq W(t,x)-W(t+\delta ,x)\leq
C(1+|x|)\delta^{\frac{1}{2}}.
 \ee
 From it we obtain immediately that $W$ is $\frac{1}{2}-$H\"{o}lder continuous in
 $t$. We will only check the second inequality in (4.30), the
 first one can be shown in a similar way. To this end we note that
 due to (4.28), for an arbitrarily small $\varepsilon>0,$
\be W(t,x)-W(t+\delta ,x) \leq I^1_\delta +I^2_\delta
+\varepsilon, \ee
 where
$$
\begin{array}{lll}
I^1_\delta & := & G^{t,x;
u^{\varepsilon},\beta(u^{\varepsilon})}_{t,t+\delta}[W(t+\delta,
X^{t,x; u^{\varepsilon},\beta(u^{\varepsilon})}_{t+\delta})]
                   -G^{t,x;u^{\varepsilon},\beta(u^{\varepsilon})}_{t,t+\delta} [W(t+\delta,x)], \\
I^2_\delta & := & G^{t,x;
u^{\varepsilon},\beta(u^{\varepsilon})}_{t,t+\delta}
[W(t+\delta,x)] -W(t+\delta ,x),
\end{array}
$$
for arbitrarily chosen $\beta\in {\cal{B}}_{t, t+\delta}$\ and
$u^{\varepsilon} \in {\cal{U}}_{t, t+\delta}$\ such that (4.28)
holds. From Lemma 2.3 and the estimate (4.12) we obtain that, for
some constant $C$ independent of the controls $u^{\varepsilon}\
\mbox{and}\ \ \beta(u^{\varepsilon})$,
$$
\begin{array}{rcl}
|I^1_\delta | &\leq& [CE(|W(t+\delta ,X^{t,x;
u^{\varepsilon},\beta(u^{\varepsilon})}_{t+\delta})
                  -W(t+\delta ,x)|^2|{{\mathcal{F}}_t})]^{\frac{1}{2}}\\
              & \leq&[CE(|X^{t,x;u^{\varepsilon},\beta(u^{\varepsilon})}_{t+\delta} -x|^2|{{\mathcal{F}}_t})]^{\frac{1}{2}},
\end{array}
$$
and since
$E[|X^{t,x;u^{\varepsilon},\beta(u^{\varepsilon})}_{t+\delta}
-x|^2|{{\mathcal{F}}_t}] \leq C(1+|x|^2) \delta $ we deduce that
$|I^1_\delta| \leq C (1+|x|)\delta^{\frac{1}{2}}$. From the
definition of $G^{t,x;
u^{\varepsilon},\beta(u^{\varepsilon})}_{t,t+\delta}[\cdot]$\ (see
(4.13)) we know that the second term $I^2_\delta $ can be written
as£º
$$
\begin{array}{llll}
I^2_\delta  & = & E[W(t+\delta ,x) +\int^{t+\delta}_t
f(s,X^{t,x;u^{\varepsilon},\beta(u^{\varepsilon})}_s,\tilde{Y}^{t,x;u^{\varepsilon},\beta(u^{\varepsilon})}_s,
\tilde{Z}^{t,x;u^{\varepsilon},\beta(u^{\varepsilon})}_s,
             u^{\varepsilon}_s, \beta_s(u^{\varepsilon}_.)) ds   \\
 &         &   -\int^{t+\delta}_t \tilde{Z}^{t,x; u^{\varepsilon},\beta(u^{\varepsilon})}_s dB_s|{{\mathcal{F}}_t}] -W(t+\delta ,x)  \\
 &    =  &  E[\int^{t+\delta}_t f(s,X^{t,x; u^{\varepsilon},\beta(u^{\varepsilon})}_s,\tilde{Y}^{t,x; u^{\varepsilon},
 \beta(u^{\varepsilon})}_s, \tilde{Z}^{t,x; u^{\varepsilon},\beta(u^{\varepsilon})}_s,
 u^{\varepsilon}_s,\beta_s(u^{\varepsilon}_.))
 ds|{{\mathcal{F}}_t}].
\end{array}
$$
With the help of the Schwartz inequality, the estimates (4.3) and
(3.4)-(i), we then have
$$
\begin{array}{lll}
|I^2_\delta | & \leq \delta^{\frac{1}{2}}
     E[\int^{t+\delta}_t |f(s,X^{t,x;u^{\varepsilon},\beta(u^{\varepsilon})}_s,
     \tilde{Y}^{t,x;u^{\varepsilon},\beta(u^{\varepsilon})}_s,\tilde{Z}^{t,x;u^{\varepsilon},\beta(u^{\varepsilon})}_s,
u^{\varepsilon}_s,\beta_s(u^{\varepsilon}_.))|^2ds|{{\mathcal{F}}_t}]^{\frac{1}{2}}  \\
& \leq\delta^{\frac{1}{2}}E[\int^{t+\delta}_t
(|f(s,X^{t,x;u^{\varepsilon},\beta(u^{\varepsilon})}_s,0,0,u^{\varepsilon}_s,\beta_s(u^{\varepsilon}_.))|+C|
\tilde{Y}^{t,x;u^{\varepsilon}, \beta(u^{\varepsilon})}_s|
+C|\tilde{Z}^{t,x;u^{\varepsilon},\beta(u^{\varepsilon})}_s|)^2ds|{{\mathcal{F}}_t}]^{\frac{1}{2}}\\
& \leq C\delta^{\frac{1}{2}}E[\int^{t+\delta}_t
(|1+|X^{t,x;u^{\varepsilon},\beta(u^{\varepsilon})}_s|+|
\tilde{Y}^{t,x;u^{\varepsilon}, \beta(u^{\varepsilon})}_s|
+|\tilde{Z}^{t,x;u^{\varepsilon},\beta(u^{\varepsilon})}_s|)^2ds|{{\mathcal{F}}_t}]^{\frac{1}{2}}\\
 & \leq C (1+|x|)\delta^{\frac{1}{2}}.
\end{array}
$$
Hence, from (4.31), $$W(t,x)-W(t+\delta ,x) \leq C
(1+|x|)\delta^{\frac{1}{2}} +\varepsilon,$$ and letting
$\varepsilon \downarrow 0$\ we get the second inequality of
(4.30). The proof is complete.\endpf

\section{\large Viscosity Solution of Isaacs' Equation: Existence Theorem }

 \hskip1cm In this section we consider the following Isaacs'
equations \be
 \left \{\begin{array}{ll}
 &\!\!\!\!\! \frac{\partial }{\partial t} W(t,x) +  H^{-}(t, x, W, DW, D^2W)=0,
 \hskip 0.5cm   (t,x)\in [0,T)\times {\mathbb{R}}^n ,  \\
 &\!\!\!\!\!  W(T,x) =\Phi (x), \hskip0.5cm   x \in {\mathbb{R}}^n,
 \end{array}\right.
\ee
and
 \be
 \left \{\begin{array}{ll}
 &\!\!\!\!\! \frac{\partial }{\partial t} U(t,x) +  H^{+}(t, x, U, DU, D^2U)=0,
 \hskip 0.5cm   (t,x)\in [0,T)\times {\mathbb{R}}^n ,  \\
 &\!\!\!\!\!  U(T,x) =\Phi (x), \hskip0.5cm   x \in {\mathbb{R}}^n,
 \end{array}\right.
\ee associated with the Hamiltonians $$ H^-(t, x, y, p, X)=
\mbox{sup}_{u \in U}\mbox{inf}_{v \in
V}\{\frac{1}{2}tr(\sigma\sigma^{T}(t, x,
 u, v)X)+ p.b(t, x, u, v)+ f(t, x, y, p.\sigma,
u, v)\}$$ and
$$ H^+(t, x, y, p, X)= \mbox{inf}_{v \in
V}\mbox{sup}_{u \in U}\{\frac{1}{2}tr(\sigma\sigma^{T}(t, x,
 u, v)X)+ p.b(t, x, u, v)+ f(t, x, y, p.\sigma,
u, v)\},$$ $\mbox{respectively, where}\ t\in [0, T],\ x\in
{\mathbb{R}}^n,\ y\in {\mathbb{R}},\ p\in {\mathbb{R}}^n\
\mbox{and}\ X\in {\mathbf{S}}^n$ $(\mbox{recall that}\
{\mathbf{S}}^n\ \mbox{denotes the}$\\ $\mbox{set of} \ n\times n \
\mbox{symmetric matrices})$. Here the functions $b, \sigma, f\
\mbox{and}\ \Phi$\ are supposed to satisfy (H4.1) and (H4.2),
respectively.

 In this section we want to prove that
the lower value function $W(t, x)$ introduced by (4.9) is the
 viscosity solution of equation (5.1), while the upper value
function $U(t, x)$ defined by (4.10) is the viscosity solution of
equation (5.2). For this we translate Peng's BSDE approach [14]
developed in the framework of stochastic control theory into that
of the stochastic differential games. Uniqueness of the viscosity
solution will be shown in the next section for the class of
continuous functions satisfying some growth assumption which is
weaker than the polynomial growth condition. We first recall the
definition of a viscosity solution of equation (5.1), similarly
for equation (5.2). The reader more interested in viscosity
solutions is referred to Crandall, Ishii and Lions [5].

\bde\mbox{ } A real-valued
continuous function $W\in C([0,T]\times {\mathbb{R}}^n )$ is called \\
  {\rm(i)} a viscosity subsolution of equation (5.1) if $W(T,x) \leq \Phi (x), \mbox{for all}\ x \in
  {\mathbb{R}}^n$, and if for all functions $\varphi \in C^3_{l, b}([0,T]\times
  {\mathbb{R}}^n)$ and $(t,x) \in [0,T) \times {\mathbb{R}}^n$ such that $W-\varphi $\ attains its
  local maximum at $(t, x)$,
     $$
     \frac{\partial \varphi}{\partial t} (t,x) +  H^{-}(t, x, \varphi, D\varphi, D^2\varphi) \geq 0;
     $$
{\rm(ii)} a viscosity supersolution of equation (5.1) if $W(T,x)
\geq \Phi (x), \mbox{for all}\ x \in
  {\mathbb{R}}^n$, and if for all functions $\varphi \in C^3_{l, b}([0,T]\times
  {\mathbb{R}}^n)$ and $(t,x) \in [0,T) \times {\mathbb{R}}^n$ such that $W-\varphi $\ attains its
  local minimum at $(t, x)$,
     $$
     \frac{\partial \varphi}{\partial t} (t,x) +  H^{-}(t, x, \varphi, D\varphi, D^2\varphi) \leq 0;
     $$
 {\rm(iii)} a viscosity solution of equation (5.1) if it is both a viscosity sub- and a supersolution of equation
     (5.1).\ede
\br \mbox{  }$C^3_{l, b}([0,T]\times {\mathbb{R}}^n)$ denotes the
set of the real-valued functions that are continuously
differentiable up to the third order and whose derivatives of
order from 1 to 3 are bounded.\er

We first prove that the lower value function $W(t,x)$ is a
viscosity solution of equation (5.1).
 \bt Under the assumptions (H4.1) and (H4.2) the lower value function $W(t,x)$ is a viscosity solution of
equation (5.1). \et

\noindent For the proof of this theorem we need four auxiliary
lemmata. To abbreviate notation we put, for some arbitrarily
chosen but fixed $\varphi \in C^3_{l, b} ([0,T] \times
{\mathbb{R}}^n)$,
 \be
\begin{array}{lll}
     F(s,x,y,z,u, v)=&\!\!\!\! \frac{\partial }{\partial s}\varphi (s,x) +
     \frac{1}{2}tr(\sigma\sigma^{T}(s,x, u, v)D^2\varphi)+ D\varphi.b(s, x, u, v) \\
        &+ f(s, x, y+\varphi (s,x), z+ D\varphi (s,x).\sigma(s,x,u, v),u, v), \\
     \end{array}
\ee $(s,x,y,z,u, v)\in [0,T] \times {\mathbb{R}}^n \times
{\mathbb{R}} \times {\mathbb{R}}^d \times U \times V,$\ and we
consider the following BSDE defined on the interval $[t,t+\delta]\
(0<\delta\leq T-t):$
    \be \left \{\begin{array}{rl}
      -dY^{1,u,v}_s =&\!\!\!\! F(s,X^{t,x;u,v}_s, Y^{1,u,v}_s,Z^{1,u,v}_s,u_s,v_s)ds
                   -Z^{1,u,v}_s dB_s, \\
     Y^{1,u,v}_{t+\delta}=&\!\!\!\!0,
     \end{array}\right.
     \ee
     where the process $X^{t,x,u,v}$\ has been introduced by equation
     $(4.1)$\ and $u(\cdot) \in {\mathcal{U}}_{t, t+\delta},\ v(\cdot) \in
     {\mathcal{V}}_{t, t+\delta}$.
\br\mbox{}It's not hard to check that $F(s,X^{t,x;u,v}_s,
y,z,u_s,v_s)$\ satisfies (A1) and (A2). Thus, due to Lemma 2.1
equation (5.4) has a unique solution. \er

      We can characterize the solution process $Y^{1,u,v}$ as follows:
    \bl\mbox{  } For every $s\in [t,t+\delta]$, we have the following
relationship:
    \be
     Y^{1,u,v}_s = G^{t,x;u,v}_{s,t+\delta} [\varphi (t+\delta ,X^{t,x;u,v}_{t+\delta})]
                -\varphi (s,X^{t,x;u,v}_s), \hskip 0.5cm
              \mbox{ P-a.s.}.\ee
 \el
\noindent \textbf{Proof}: We recall that $G^{t,x;u,v}_{s,t+\delta}
[\varphi (t+\delta, X^{t,x;u,v}_{t+\delta})]$ is defined with the
help of the solution of the BSDE
     $$
     \left \{\begin{array}{rl}
     -dY^{u,v}_s =\!\!\! & f(s,X^{t,x;u,v}_s, Y^{u,v}_s,Z^{u,v}_s,u_s,v_s)ds
      -Z^{u,v}_s dB_s , \hskip 0.2cm s\in [t,t+\delta], \\
     Y^{u,v}_{t+\delta}=\!\!\!& \varphi (t+\delta
     ,X^{t,x;u,v}_{t+\delta}),
     \end{array}\right.
     $$
by the following formula:
     \be
     G^{t,x;u,v}_{s,t+\delta} [\varphi (t+\delta ,X^{t,x;u,v}_{t+\delta})] =Y^{u,v}_s, \hskip 0.5cm
     s\in [t,t+\delta]  \ee
(see (4.13)). Therefore we only need to prove that
$Y^{u,v}_s-\varphi (s,X^{t,x;u,v}_s)\equiv Y^{1,u,v}_s.$\ This
result can be obtained easily by applying It$\hat{o}$'s formula to
$\varphi (s,X^{t,x;u,v}_s)$. Indeed, we get that the stochastic
differentials of $Y^{u,v}_s -\varphi (s,X^{t,x;u,v}_s)$\ and
$Y^{1,u,v}_s$\ coincide, while at the terminal time $t+\delta$,
$Y^{u,v}_{t+\delta} - \varphi (t+\delta ,X^{t,x;u,v}_{t+\delta})
=0 = Y^{1,u,v}_{t+\delta}.$
     So the proof is complete.\endpf \vskip 0.3cm

Now we consider the following simple BSDE in which the driving
process $X^{t,x;u,v}$ is replaced by its deterministic initial
value $x$: \be
     \left \{\begin{array}{rl}
     -dY^{2,u,v}_s=&\!\!\!  F(s,x,Y^{2,u,v}_s ,Z^{2,u,v}_s,u_s, v_s)ds -Z^{2,u,v}_s dB_s,
         \\
     Y^{2,u,v}_{t+\delta}=&\!\!\! 0,
         \hskip 0.3cm s\in [t,t+\delta],
     \end{array}\right.
 \ee
where $u(\cdot) \in {\mathcal{U}}_{t, t+\delta},\ v(\cdot) \in
     {\mathcal{V}}_{t, t+\delta}$.
The following Lemma will allow us to neglect the difference
$|Y^{1,u,v}_t-Y^{2,u,v}_t|$ for sufficiently small $\delta >0$.

\bl For every $u \in {\mathcal{U}}_{t, t+\delta},\ v \in
{\mathcal{V}}_{t, t+\delta},$\ we have
 \be |Y^{1,u,v}_t-Y^{2,u,v}_t| \leq C\delta^{\frac{3}{2}},\ \ \mbox{P-a.s.},
       \ee
 where C is independent
 of the control processes $u$\ and $v$.
\el
 \noindent \textbf{Proof}: From (4.3) we have for all $p\geq 2$ the existence
 of some $C_{p}\in {\mathbb{R}}_+$\ such that
    $$
    E [\sup\limits_{t\leq s \leq T} |X^{t,x;u,v}_s|^p|{{\mathcal{F}}_t}]\leq
    C_{p}(1+|x|^p),\ \
    \mbox{P-a.s., \ uniformly in}\ u \in {\mathcal{U}}_{t, t+\delta}, v \in
{\mathcal{V}}_{t, t+\delta}.$$
  This combined with the estimate
     $$
      \begin{array}{lll}
     E [\sup \limits_{t\leq s \leq t+\delta} |X^{t,x;u,v}_s -x|^p |{{\mathcal{F}}_t}]& \leq &
     2^{p-1}E [\sup \limits_{t\leq s \leq t+\delta} |\int^s_t b(r,X^{t,x;u,v}_r,u_r, v_r)dr|^p|{{\mathcal{F}}_t}]  \\
     & &+  2^{p-1} E[ \sup \limits_{t\leq s \leq t+\delta} |\int^s_t \sigma (r,X^{t,x;u,v}_r,u_r, v_r)dB_r|^p|{{\mathcal{F}}_t}] \\
     \end{array}
     $$
    yields
      \be
       E [\sup \limits_{t\leq s \leq t+\delta} |X^{t,x;u,v}_s -x|^p|{{\mathcal{F}}_t}] \leq
       C_p\delta^{\frac{p}{2}},\ \
    \mbox{P-a.s., \ uniformly in}\ u \in {\mathcal{U}}_{t, t+\delta}, v \in
{\mathcal{V}}_{t, t+\delta}.
      \ee
   We now apply Lemma 2.3 combined with (5.9) to equations (5.4) and (5.7). For this we
set in Lemma 2.3:
       $$\xi_1 =\xi_2 =0,\ g(s,y,z) =F(s,X^{t,x,u,v}_s,y,z,u_s,v_s),$$
       $$\varphi_1(s)=0,\ \varphi_2(s)=F(s,x,Y_s^{2, u, v},Z_s^{2, u, v},u_s,v_s)-F(s,X^{t,x,u,v}_s,Y_s^{2, u, v},Z_s^{2, u, v},u_s,v_s).
       $$
       Obviously, the function $g$\ is Lipschitz with respect to $(y,z)$, and
       $|\varphi_2(s)|\leq C(1+|x|^2)(|X^{t,x;u,v}_s -x|+|X^{t,x;u,v}_s -x|^3),$\
       for $s\in [t, t+\delta], (t, x)\in [0, T)\times {\mathbb{R}}^n$, $u\in {\mathcal{U}}_{t, t+\delta}, v \in {\mathcal{V}}_{t, t+\delta}.
       $
       Thus, with the notation $\rho_0 (r) =(1+|x|^2)(r+r^3),\ r\geq 0, $\ we have
      $$
       \begin{array}{lll}
       &  & E[\int^{t+\delta}_t (|Y^{1,u,v}_s -Y^{2,u,v}_s|^2 +|Z^{1,u,v}_s -Z^{2,u,v}_s|^2)ds|{\mathcal{F}}_{t}]  \\
       & \leq & CE[\int^{t+\delta}_t \rho^2_0 (|X^{t,x,u,v}_s -x|) ds|{\mathcal{F}}_{t}]\\
       & \leq &C \delta E[\sup \limits_{t\leq s \leq t+\delta}\rho^2_0 (|X^{t,x,u,v}_s
       -x|)|{\mathcal{F}}_{t}]\\
       &\leq & C\delta^2.
       \end{array}
     $$
Therefore,
       $$
       \begin{array}{llll}
      && |Y^{1,u,v}_t -Y^{2,u,v}_t|  =|E[(Y^{1,u,v}_t -Y^{2,u,v}_t )|{\mathcal{F}}_{t}]|  \\
       & = & |E[\int^{t+\delta}_t (F(s,X^{t,x,u,v}_s,Y^{1,u,v}_s,Z^{1,u,v}_s,u_s, v_s)
                -F(s,x,Y^{2,u,v}_s ,Z^{2,u,v}_s,u_s, v_s)) ds|{\mathcal{F}}_{t}]|   \\
       & \leq & CE [\int^{t+\delta}_t [\rho_0 (|X^{t,x,u,v}_s -x|) +|Y^{1,u,v}_s -Y^{2,u,v}_s|
       +|Z^{1,u,v}_s -Z^{2,u,v}_s|] ds|{\mathcal{F}}_{t}]  \\
       & \leq & CE[\int^{t+\delta}_t\rho_0 (|X^{t,x,u,v}_s -x|)ds|{\mathcal{F}}_{t}] + C\delta^{\frac{1}{2}} \{ E[\int^{t+\delta}_t
            |Y^{1,u,v}_s
            -Y^{2,u,v}_s|^2|{\mathcal{F}}_{t}]^{\frac{1}{2}}\\
         &&   + E[\int^{t+\delta}_t|Z^{1,u,v}_s -Z^{2,u,v}_s|^2ds|{\mathcal{F}}_{t}]^{\frac{1}{2}}\}  \\
       & \leq & C\delta^{\frac{3}{2}}.
       \end{array}
       $$
   Thus, the proof is
   complete.\endpf
     \vskip 0.3cm
\bl Let $Y_0 (\cdot)$ be the solution of the following ordinary
differential equation:
     \be
     \left \{\begin{array}{lll}
     - {\dot{Y}}_0 (s)  & =  & F_0(s,x,Y_0 (s),0),\ \ s\in [t,t+\delta], \\
      Y_0 (t+\delta )   & =  & 0,
     \end{array}\right.
    \ee
  where the function $F_0$ is defined by
    \be
     F_0(s,x,y,z)=\mbox{sup}_{u\in U}\mbox{inf}_{v \in
V} F(s,x,y,z,u,v). \ee Then, P-a.s.,\be\mbox{esssup}_{u \in
{\mathcal{U}}_{t, t+\delta}}\mbox{essinf}_{v \in {\mathcal{V}}_{t,
t+\delta}} Y^{2,u,v}_t =Y_0 (t). \ee
  \el
 \noindent \textbf{Proof}: Obviously,
 $F_0(s,x,y,z)$\ is Lipschitz in $(y, z)$, uniformly with respect to $(s, x).$\ This guarantees existence and uniqueness for equation
 (5.10). We first introduce the function \be
     F_1(s,x,y,z,u)=\mbox{inf}_{v\in V} F(s,x,y,z,u,v),\ (s,x,y,z,u)\in [0, T]\times {\mathbb{R}}^n \times {\mathbb{R}}\times {\mathbb{R}}^d \times U, \ee
and consider the BSDE \be
     \left \{\begin{array}{lll}
     - d{Y}^{3,u} (s)  & =  & F_1(s,x,{Y}^{3,u}(s),{Z}^{3,u}(s),u_s)ds-{Z}^{3,u}(s)dB_s, \\
      {Y}^{3,u} (t+\delta )   & =  & 0,\ \ \ s\in [t,t+\delta], \ \
       \end{array}\right.
    \ee
$\mbox{for}\ u\in {\mathcal{U}}_{t, t+\delta}$. We notice that
since $F_1(s,x,y,z,u_s)$\ is Lipschitz in $(y, z)$, for every
$u\in {\mathcal{U}}_{t, t+\delta},$\ there exists a unique
solution $({Y}^{3,u}, {Z}^{3,u})$\ to the BSDE (5.14). Moreover,
$${Y}^{3,u} (t)= \mbox{essinf}_{v(\cdot) \in {\mathcal{V}}_{t, t+\delta}}
Y^{2,u,v}_t,\ \mbox{P-a.s.},\ \mbox{for any}\ u\in
      {\mathcal{U}}_{t, t+\delta}.$$
Indeed, from the definition of $F_1$\ and Lemma 2.2 (comparison
theorem) we have
$${Y}^{3,u} (t)\leq \mbox{essinf}_{v(\cdot) \in {\mathcal{V}}_{t, t+\delta}}
Y^{2,u,v}_t,\ \mbox{P-a.s.},\ \mbox{for all}\ u\in
      {\mathcal{U}}_{t, t+\delta}.$$
On the other hand, there exists a measurable function $v^3: [t,
T]\times
{\mathbb{R}}^n\times{\mathbb{R}}\times{\mathbb{R}}^d\times
U\rightarrow V$\ such that
$$F_1(s,x,y,z,u)= F(s,x,y,z,u,v^3(s,x,y,z,u)),\ \mbox{for any}\ s, x, y, z, u.$$
 We then put
$$\widetilde{v}_s^3:=
     v^{3} (s, x, {Y}^{3,u}_s, {Z}^{3,u}_s, u_s ), \  s\in [t,
     t+\delta],
$$
and we observe that $\widetilde{v}^3 \in {\mathcal{V}}_{t,
t+\delta},$\ and
$$F_1(s,x,{Y}^{3,u}_s, {Z}^{3,u}_s, u_s)= F(s,x,{Y}^{3,u}_s, {Z}^{3,u}_s, u _s, \widetilde{v}_s^3),\ s\in [t, t+\delta].$$
Consequently, from the uniqueness of the solution of the BSDE it
follows that $({Y}^{3,u}, {Z}^{3,u})=({Y}^{2,u,\widetilde{v}^3},
{Z}^{2,u, \widetilde{v}^3})$\ and, in particular,
${Y}^{3,u}_t={Y}^{2,u,\widetilde{v}^3}_t,\ \mbox{P-a.s.},
\mbox{for any}\ u\in {\mathcal{U}}_{t,t+\delta}.$\ This proves
that
$${Y}^{3,u} (t)= \mbox{essinf}_{v\in {\mathcal{V}}_{t, t+\delta}}
Y^{2,u,v}_t,\ \mbox{P-a.s.},\ \mbox{for all}\ u\in
      {\mathcal{U}}_{t, t+\delta}.$$
 Finally, since $F_0(s,x,y,z)=\mbox{sup}_{u \in
U}F_1(s,x,y,z,u),$ \ an argument similar to that developed above
yields
$$Y_0 (t)=\mbox{esssup}_{u \in
{\mathcal{U}}_{t, t+\delta}} Y^{3,u}(t)(=\mbox{esssup}_{u \in
{\mathcal{U}}_{t, t+\delta}}\mbox{essinf}_{v \in {\mathcal{V}}_{t,
t+\delta}} Y^{2,u,v}_t), \ \mbox{P-a.s.}.$$ It uses the fact that
equation (5.10) can be considered as a BSDE with solution $(Y_s,
Z_s)= (Y_0(s), 0).$
 The proof is complete.\endpf
\bl For every $u \in {\mathcal{U}}_{t, t+\delta}, v \in
{\mathcal{V}}_{t, t+\delta},$\ we have\be
E[\int_t^{t+\delta}|Y^{2,u,v}_s|ds|{\mathcal{F}}_{t}]+E[\int_t^{t+\delta}|Z^{2,u,v}_s|ds|{\mathcal{F}}_{t}]\leq
C\delta^{\frac{3}{2}},\ \mbox{P-a.s.}, \ee  where the constant
$C$\ is independent of the controls $u,\ v$.\el
 \noindent \textbf{Proof}: Since $F(s, x, \cdot, \cdot, u, v)$\ has a linear
 growth in $(y, z)$, uniformly in $(u, v)$, we get from Lemma 2.3 that, for some constant $C$\ independent of $\delta$\
 and the control processes $u,\ v,$
 $$|Y^{2,u,v}_s|^2\leq C\delta,\ E[\int_s^{t+\delta}|Z^{2,u,v}_r|^2dr|{\cal{F}}_s]\leq
C\delta,\ s\in[t, t+\delta].$$
On the other hand, from equation
(5.7),
 $$\begin{array}{ll}
 |Y^{2,u,v}_s|&\leq E[\int_s^{t+\delta}|F(r,x,{Y}^{2, u, v}_r, {Z}^{2, u, v}_r, u _r, v_r)|dr|{\cal{F}}_s]\\
 &\leq CE[\int_s^{t+\delta}(1+ |x|^2+|{Y}^{2, u, v}_r|+ |{Z}^{2, u, v}_r|)dr|{\cal{F}}_s]\\
 &\leq C\delta+C\sqrt{\delta}(E[\int_s^{t+\delta}|Z^{2,u,v}_r|^2dr|{\cal{F}}_s])^{\frac{1}{2}}\leq C\delta,\ \mbox{P-a.s.},s\in[t,
 t+\delta],
 \end{array}$$
and, since
$$\int_t^{t+\delta}{Z}^{2, u, v}_sdB_s=\int_t^{t+\delta}F(s,x,{Y}^{2, u, v}_s, {Z}^{2, u, v}_s, u _s, v_s)ds-Y^{2,u,v}_t,$$
we can get $E[\int_t^{t+\delta}|Z^{2,u,v}_s|^2ds|{\cal{F}}_t]\leq
C\delta^2.$\ Finally,
$$\begin{array}{ll}E[\int_t^{t+\delta}|Y^{2,u,v}_s|ds|{\cal{F}}_t]+E[\int_t^{t+\delta}|Z^{2,u,v}_s|ds|{\cal{F}}_t]&\leq
C\delta^2+\delta^{\frac{1}{2}}\{E[\int_t^{t+\delta}|Z^{2,u,v}_s|^2ds|{\cal{F}}_t]\}^{\frac{1}{2}}\\\leq
C\delta^{\frac{3}{2}}.\end{array}$$ The proof is complete.\endpf
\vskip0.2cm
 Now we are able to give the proof of Theorem 5.1:

\noindent \textbf{Proof}: (1) Obviously, $W(T,x)=\Phi (x),\ x\in
{\mathbb{R}}^n $. Let us show in a first step that $W$ is a
viscosity supersolution. For this we suppose that $\varphi \in
C^3_{l,b} ([0,T] \times {\mathbb{R}}^n)$,\ and $(t,x)\in [0,
T)\times {\mathbb{R}}^n$\ are such that $W-\varphi$\ attains its
minimum at $(t,x).$\ Notice that we can replace the condition of a
local minimum by that of a global one in the definition of the
viscosity supersolution since $W$ is continuous and of at most
linear growth. Without loss of generality we may also suppose that
$\varphi (t,x)=W(t,x)$.\ Then, due to the DPP (see Theorem 4.2),
     $$
     \varphi (t,x) =W(t,x) =\mbox{essinf}_{\beta \in {\mathcal{B}}_{t, t+\delta}}\mbox{esssup}_{u \in
{\mathcal{U}}_{t, t+\delta}}G^{t,x;u,\beta(u)}_{t,t+\delta}
[W(t+\delta, X^{t,x;u,\beta(u)}_{t+\delta})],\ 0\leq\delta\leq
T-t,
     $$
 and from $W\geq \varphi$\ and the monotonicity property of $G^{t,x;u,\beta(u)}_{t,t+\delta}[\cdot]$\ (see Lemma 2.2)  we obtain
     $$
     \mbox{essinf}_{\beta \in {\mathcal{B}}_{t, t+\delta}}\mbox{esssup}_{u \in
{\mathcal{U}}_{t, t+\delta}} \{G^{t,x;u,\beta(u)}_{t,t+\delta}
[\varphi(t+\delta, X^{t,x;u,\beta(u)}_{t+\delta})] -\varphi
(t,x)\}\leq 0,\ \mbox{P-a.s.}.
     $$
 Thus, from Lemma 5.1,
    $$
      \mbox{essinf}_{\beta \in {\mathcal{B}}_{t, t+\delta}}\mbox{esssup}_{u \in
{\mathcal{U}}_{t, t+\delta}} Y^{1,u,\beta(u)}_t \leq 0,\
\mbox{P-a.s.},
    $$
and further, from Lemma 5.2 we have
     $$
     \mbox{essinf}_{\beta \in {\mathcal{B}}_{t, t+\delta}}\mbox{esssup}_{u \in
{\mathcal{U}}_{t, t+\delta}} Y^{2,u,\beta(u)}_t \leq
C\delta^{\frac{3}{2}},\ \mbox{P-a.s.}.
     $$
     Consequently, since $\mbox{essinf}_{v \in
{\mathcal{V}}_{t, t+\delta}} Y^{2,u,v}_t\leq Y^{2,u,\beta(u)}_t, \
\beta\in{\mathcal{B}}_{t, t+\delta} $, we get
$$
\mbox{esssup}_{u \in {\mathcal{U}}_{t, t+\delta}}\mbox{essinf}_{v
\in {\mathcal{V}}_{t, t+\delta}}
Y^{2,u,v}_t\leq\mbox{essinf}_{\beta \in {\mathcal{B}}_{t,
t+\delta}}\mbox{esssup}_{u \in {\mathcal{U}}_{t, t+\delta}}
Y^{2,u,\beta(u)}_t \leq C\delta^{\frac{3}{2}},\ \mbox{P-a.s.},
$$
and Lemma 5.3 implies
     $$
     Y_0 (t) \leq C\delta^{\frac{3}{2}},\ \mbox{P-a.s.},
     $$
 where $Y_0$ is the unique solution of equation (5.10). It then
 follows easily that
     $$
      \mbox{sup}_{u\in U} \mbox{inf}_{v \in
V}F(t,x,0,0,u, v)=F_0 (t,x,0,0)\leq 0,
     $$
and from the definition of $F$\ we see that $W$ is a viscosity
 supersolution of equation (5.1).\\

(2) The second step is devoted to the proof that $W$ is a
viscosity subsolution. For this we suppose that $\varphi \in
C^3_{l,b} ([0,T] \times {\mathbb{R}}^n)$\ and $(t,x)\in [0,
T)\times {\mathbb{R}}^n$\ are such that $W-\varphi$\ attains its
maximum at $(t,x)$. Without loss of generality we suppose again
$\varphi (t,x)=W(t,x)$. We must prove that
$$      \mbox{sup}_{u\in U} \mbox{inf}_{v \in
V}F(t,x,0,0,u, v)=F_0 (t,x,0,0)\geq 0.
     $$
Let us suppose that this is not true. Then there exists some
$\theta>0$ such that \be F_0 (t,x,0,0)=\mbox{sup}_{u\in U}
\mbox{inf}_{v \in V}F(t,x,0,0,u, v)\leq-\theta<0,\ee and we can
find a measurable function $\psi: U\rightarrow V$ such that
$$ F(t,x,0,0,u,\psi(u))\leq -\frac{3}{4}\theta,\ \mbox{for all}\
u\in U.$$ Moreover, since $F(\cdot,x,0,0,\cdot,\cdot)$\ is
uniformly continuous on $[0, T]\times U\times V$\ there exists
some $T-t\geq R>0$\ such that \be F(s,x,0,0,u,\psi(u))\leq
-\frac{1}{2}\theta,\ \mbox{for all}\ u\in U \mbox{and}\ |s-t|\leq
R. \ee On the other hand, due to the DPP (see Theorem 4.1), for
every $\delta\in(0, R]$,
     $$
     \varphi (t,x) =W(t,x) =\mbox{essinf}_{\beta \in {\mathcal{B}}_{t, t+\delta}}\mbox{esssup}_{u \in
{\mathcal{U}}_{t, t+\delta}}G^{t,x;u,\beta(u)}_{t,t+\delta}
[W(t+\delta, X^{t,x;u,\beta(u)}_{t+\delta})],
     $$
  and from $W\leq\varphi$\ and the monotonicity property of $G^{t,x;u,\beta(u)}_{t,t+\delta}[\cdot]$\ (see Lemma 2.2)\ we obtain
     $$
     \mbox{essinf}_{\beta \in {\mathcal{B}}_{t, t+\delta}}\mbox{esssup}_{u \in
{\mathcal{U}}_{t, t+\delta}} \{G^{t,x;u,\beta(u)}_{t,t+\delta}
[\varphi(t+\delta, X^{t,x;u,\beta(u)}_{t+\delta})] -\varphi
(t,x)\}\geq 0,\ \mbox{P-a.s.}.
     $$
 Thus, from Lemma 5.1,
    $$
      \mbox{essinf}_{\beta \in {\mathcal{B}}_{t, t+\delta}}\mbox{esssup}_{u \in
{\mathcal{U}}_{t, t+\delta}} Y^{1,u,\beta(u)}_t \geq 0,\
\mbox{P-a.s.},
    $$
and, in particular,
    $$
     \mbox{esssup}_{u \in
{\mathcal{U}}_{t, t+\delta}} Y^{1,u,\psi(u)}_t \geq 0,\
\mbox{P-a.s.}.
    $$
 Here, by putting $\psi_s(u)(\omega)=\psi(u_s(\omega)),\ (s,\omega)\in[t,T]\times
 \Omega$, we identify $\psi$\ as an element of ${\cal{B}}_{t, t+\delta}$.
Given an arbitrarily $\varepsilon>0$ we can choose
$u^\varepsilon\in {\mathcal{U}}_{t, t+\delta}$\ such that
$Y^{1,u^\varepsilon,\psi(u^\varepsilon)}_t\geq
-\varepsilon\delta.$\ For this the argument developed in the proof
of the lemmata can be used. From Lemma 5.2 we further have
     \be
   Y^{2,u^\varepsilon,\psi(u^\varepsilon)}_t\geq -C\delta^{\frac{3}{2}}- \varepsilon\delta,\ \mbox{P-a.s.}.
     \ee
  Taking into account that $Y^{2,u^\varepsilon,\psi(u^\varepsilon)}_t
  =E[\int_t^{t+\delta}F(s, x, Y^{2,u^\varepsilon,\psi(u^\varepsilon)}_s,Z^{2,u^\varepsilon,\psi(u^\varepsilon)}_s,
  u^\varepsilon_s, \psi_s(u^\varepsilon_.))ds|{\mathcal{F}}_{t}]$\ we get
from the Lipschitz property of $F$ in $(y, z)$, (5.17) and Lemma
5.4 that \be\begin{array}{ll}
Y^{2,u^\varepsilon,\psi(u^\varepsilon)}_t&\leq
E[\int_t^{t+\delta}(C|Y^{2,u^\varepsilon,\psi(u^\varepsilon)}_s|+C|Z^{2,u^\varepsilon,\psi(u^\varepsilon)}_s|+F(s,
x,0,0,u^\varepsilon_s, \psi_s(u^\varepsilon_.)))ds|{\mathcal{F}}_{t}]\\
&\leq C\delta^{\frac{3}{2}}-\frac{1}{2}\theta\delta,\
\mbox{P-a.s.}.\end{array}\ee From (5.18) and (5.19),
$-C\delta^{\frac{1}{2}}- \varepsilon\leq
C\delta^{\frac{1}{2}}-\frac{1}{2}\theta,\ \mbox{P-a.s.}.$\ Letting
$\delta\downarrow0$,\ and then $\varepsilon\downarrow0$\ we deduce
$\theta\leq0$\ which induces a contradiction.
 Therefore,
     $$
     F_0 (t,x,0,0) = \mbox{sup}_{u\in U} \mbox{inf}_{v \in
V}F(t,x,0,0,u, v) \geq 0, $$
 and from the definition of $F$, we know that $W$ is a viscosity
 subsolution of equation (5.1). Finally, the results from the first and the second step prove that $W$ is a viscosity
 solution of equation (5.1).

\br Similarly, we can prove that $U$\ is a viscosity
 solution of equation (5.2). \er

\section{\large Viscosity Solution of Isaacs' Equation: Uniqueness Theorem }
The objective of this section is to study the uniqueness of the
viscosity solution of Isaacs' equation (5.1), \be
 \left \{\begin{array}{ll}
 &\!\!\!\!\! \frac{\partial }{\partial t} \omega(t,x) +  H^{-}(t, x, \omega, D\omega, D^2\omega)=0,
 \hskip 0.5cm  (t,x)\in [0,T)\times {\mathbb{R}}^n ,  \\
 &\!\!\!\!\!  \omega(T,x) =\Phi (x), \hskip0.5cm  x \in
 {\mathbb{R}}^n.
 \end{array}\right.
\ee Recall that
$$H^{-}(t, x, y, p, X)=\mbox{sup}_{u \in
U} \mbox{inf}_{v \in V}\{\frac{1}{2}tr(\sigma\sigma^{T}(t, x,
 u, v)X)+ p.b(t, x, u, v)+ f(t, x, y, p.\sigma,
u, v)\},$$ $  t\in [0, T],\ x\in {\mathbb{R}}^n,\ y\in
{\mathbb{R}},\ p\in {\mathbb{R}}^n,\ X\in {\mathbf{S}}^n$.
 The functions $b, \sigma, f\ \mbox{and}\ \Phi$\ are still supposed to satisfy (H4.1) and (H4.2), respectively.

 We will prove the uniqueness for equation (6.1) in the following
 space of continuous functions

 $\Theta=\{ \varphi\in C([0, T]\times {\mathbb{R}}^n): \exists\ \widetilde{A}>0\ \mbox{such
 that}$ \vskip 0.1cm
 $\mbox{ }\hskip2cm \lim_{|x|\rightarrow \infty}\varphi(t, x)\exp\{-\widetilde{A}[\log((|x|^2+1)^{\frac{1}{2}})]^2\}=0,\
 \mbox{uniformly in}\ t\in [0, T]\}.$  \vskip 0.1cm
\noindent This space of continuous functions endowed with a growth
condition which is slightly weaker than the assumption of
polynomial growth but more restrictive than that of exponential
growth. This growth condition was introduced by Barles, Buckdahn,
Pardoux~\cite{BBE} to prove the uniqueness of the viscosity
solution of an integro-partial differential equation associated
with a decoupled FBSDE with jumps. It was shown in~\cite{BBE} that
this kind of growth condition is optimal for the uniqueness and
can, in general, not be weakened. We adapt the ideas developed
in~\cite{BBE} to Isaacs' equation (6.1) to prove the uniqueness of
the viscosity solution in $\Theta$. Since the proof of the
uniqueness in $\Theta$\ for equation (5.2) is the same we will
restrict ourselves only on that of (6.1). Before stating the main
result of this section, let us begin with two auxiliary lemmata.
Denoting by $K$\ a Lipschitz constant of $f(t,x,.,.)$, that is
uniformly in $(t, x),$\ we have the following

\bl\mbox{  } Let $u_1 \in \Theta$\ be a viscosity subsolution and
$u_2 \in \Theta$\ be a viscosity supersolution of equation (6.1).
Then the function $\omega:= u_1-u_2$\ is a viscosity subsolution
of the equation
    \be
    \left \{
     \begin{array}{lll}
     &\!\!\!\!\!\frac{\partial }{\partial t} \omega(t,x) + \mbox{sup}_{u \in
U, v \in V}\{ \frac{1}{2}tr(\sigma\sigma^{T}(t, x,
 u, v)D^2\omega)+ D\omega.b(t, x, u, v)+ K|\omega| +\\
 &\!\!\!\!\!\mbox{ }\hskip1cm K|D\omega.\sigma(t, x, u, v)| \}= 0, \ \hskip2cm  (t, x)\in [0, T)\times
 {\mathbb{R}}^n,\\
&\!\!\!\!\!\omega(T,x) =0,\ \hskip1cm  x \in {\mathbb{R}}^n.
     \end{array}\right.
   \ee
  \el
\noindent The proof of this lemma follows directly that of Lemma
3.7 in~\cite{BBE}, it is even simpler because contrary to Lemma
3.7 in~\cite{BBE} we don't have any integral part here in equation
(6.1). In analogy to~\cite{BBE} we also have

 \bl\mbox{  }For any
$\widetilde{A}>0,$\ there exists $C_1>0$\ such that the function
$$\chi(t,x)=\exp[(C_1(T-t)+\widetilde{A})\psi(x)],$$
with
$$\psi(x)=[\log((|x|^2+1)^{\frac{1}{2}})+1]^2,\ x\in {\mathbb{R}}^n,$$
satisfies
    \be
     \begin{array}{lll}
     &\frac{\partial }{\partial t}\chi(t,x) + \mbox{sup}_{u \in
U, v \in V}\{ \frac{1}{2}tr(\sigma\sigma^{T}(t, x,
 u, v)D^2\chi)+ D\chi.b(t, x, u, v)+ K\chi(t,x) +\\
 & K|D\chi(t,x).\sigma(t, x, u, v)| \}< 0 \ \  \mbox{in}\ [t_1, T]\times
 {\mathbb{R}}^n,\ \mbox{where}\ \ t_1=T-\frac{\widetilde{A}}{C_1}.
     \end{array}
    \ee
   \el
\noindent {\bf Proof.} By direct calculus we first deduce the
following estimates for the first and second derivatives of
$\psi$:
$$ |D\psi(x)|\leq \frac{2[\psi(x)]^{\frac{1}{2}}}{(|x|^2+1)^{\frac{1}{2}}}\leq
4,\ \ \  |D^2\psi(x)|\leq
\frac{C(1+[\psi(x)]^{\frac{1}{2}})}{|x|^2+1},\ \ \ x\in
{\mathbb{R}}^n.$$ These estimates imply that, if $t\in [t_1, T],$
$$
\begin{array}{lll}
     |D\chi(t,x)|&\leq (C_1(T-t)+\widetilde{A})\chi(t,x)|D\psi(x)|\\
 & \leq
 C\chi(t,x)\frac{[\psi(x)]^{\frac{1}{2}}}{(|x|^2+1)^{\frac{1}{2}}},
     \end{array}
$$
and, similarly
$$ |D^2\chi(t,x)| \leq C\chi(t,x)\frac{\psi(x)}{|x|^2+1}.
$$
We should notice that the above estimates do not depend on $C_1$\
because of the definition of $t_1$. In virtue with the above
estimates we have
$$  \begin{array}{lll}
     &\frac{\partial }{\partial t}\chi(t,x) + \mbox{sup}_{u \in
U, v \in V}\{ \frac{1}{2}tr(\sigma\sigma^{T}(t, x,
 u, v)D^2\chi)+ D\chi.b(t, x, u, v)+ K\chi(t,x) +\\
 &\mbox{}\hskip1cm K|D\chi(t,x).\sigma(t, x, u, v)| \}\\
 &\leq -\chi(t,x)\{C_1\psi(x)-C\psi(x)-C[\psi(x)]^{\frac{1}{2}}-K\}\\
 &< -\chi(t,x)\{C_1-[2C+K]\}\psi(x)< 0,\ \mbox{if}\ C_1>2C+K \ \mbox{large
 enough}.
     \end{array}
   $$
\endpf

\noindent Now we can prove the uniqueness theorem.\\

\bt\mbox{ } We assume that (H4.1), (H4.2) hold. Let $u_1$ (resp.,
$u_2$) $\in \Theta$
  be a viscosity subsolution (resp., supersolution) of equation
  (6.1). Then we have
     \be
     u_1 (t,x) \leq u_2 (t,x) , \hskip 0.5cm \mbox{for all}\ \ (t,x) \in [0,T] \times {\mathbb{R}}^n .
    \ee
\et \noindent {\bf Proof.} Let us put $\omega:= u_1-u_2$. Then we
have, for some $\widetilde{A}>0$,
$$\lim_{|x|\rightarrow \infty}\omega(t, x)e^{-\widetilde{A}[\log((|x|^2+1)^{\frac{1}{2}})]^2}=0,$$
uniformly with respect to $t\in [0, T]$. This implies, in
particular, that for any $\alpha >0$, $\omega(t, x)-\alpha\chi(t,
x)$\ is bounded from above in $[t_1, T]\times {\mathbb{R}}^n,$\
and that
$$ M:=\max_{[t_1, T]\times {\mathbb{R}}^n}(\omega-\alpha\chi)(t, x)e^{-K(T-t)}$$
is achieved at some point $(t_0, x_0)\in [t_1, T]\times
{\mathbb{R}}^n$\ (depending on $\alpha$). We now have to distinguish between two cases.\\
For the first case we suppose that: $\omega(t_0, x_0)\leq 0$, for any $\alpha>0$.\\
Then, obviously $M\leq 0$ and $u_1(t, x)-u_2(t, x)\leq
\alpha\chi(t, x)$\ in $[t_1, T]\times {\mathbb{R}}^n$.
Consequently, letting $\alpha$\ tend to zero we obtain
$$u_1(t, x)\leq u_2(t, x),\ \ \mbox{for all}\ (t, x)\in [t_1, T]\times {\mathbb{R}}^n. $$
For the second case we assume that there exists some $\alpha>0$\ such that $\omega(t_0, x_0)> 0$.\\
We notice that $\omega(t, x)-\alpha\chi(t,x)\leq (\omega(t_0,
x_0)-\alpha\chi(t_0, x_0))e^{-K(t-t_0)}\ \ \mbox{in}\ \ [t_1,
T]\times {\mathbb{R}}^n.$\ Then, putting
$$\varphi(t, x)=\alpha\chi(t, x)+(\omega-\alpha\chi)(t_0, x_0)e^{-K(t-t_0)}$$
we get $\omega-\varphi\leq 0=(\omega-\varphi)(t_0, x_0)\
\mbox{in}\ \ [t_1, T]\times {\mathbb{R}}^n.$\ Consequently, since
$\omega$\ is a viscosity subsolution of (6.2) from Lemma 6.1  we
have
 $$
     \begin{array}{lll}
     &\frac{\partial }{\partial t}\varphi(t_0, x_0) + \mbox{sup}_{u \in
U, v \in V}\{ \frac{1}{2}tr(\sigma\sigma^{T}(t_0, x_0,
 u, v)D^2\varphi(t_0, x_0))+ D\varphi(t_0, x_0).b(t_0, x_0, u, v)+\\
  & K|\varphi(t_0, x_0)| +
 K|D\varphi(t_0, x_0).\sigma(t_0, x_0, u, v)| \}\geq 0.
     \end{array}
$$
Moreover, due to our assumption that $\omega(t_0, x_0)>0$\ and
since $\omega(t_0, x_0)=\varphi(t_0, x_0)$\ we can replace
$K|\varphi(t_0, x_0)|$\ by $K\varphi(t_0, x_0)$\ in the above
formula. Then, from the definition of $\varphi$\ and Lemma 6.2,
$$
     \begin{array}{lll}
     &0\leq \alpha\{\frac{\partial \chi}{\partial t} (t_0, x_0) + \mbox{sup}_{u \in
U, v \in V}\{ \frac{1}{2}tr(\sigma\sigma^{T}(t_0, x_0,
 u, v)D^2\chi(t_0, x_0))+ D\chi(t_0, x_0).b(t_0, x_0, u, v)+ \\
 & K\chi(t_0, x_0) + K|D\chi(t_0, x_0).\sigma(t_0, x_0, u, v)| \}\}< 0
     \end{array}
   $$
which is a contradiction. Finally, by applying successively the
same argument on the interval $[t_2, t_1]$\ with
$t_2=(t_1-\frac{\widetilde{A}}{C_1})^{+},$\ and then, if $t_2>0,$\
on $[t_3, t_2]$\ with $t_3=(t_2-\frac{\widetilde{A}}{C_1})^{+},$\
etc. We get
$$  u_1 (t,x) \leq u_2 (t,x) , \hskip 0.5cm (t,x) \in [0,T] \times {\mathbb{R}}^n .$$
Thus, the proof is complete.\endpf

   \br\mbox{  } Obviously, since the lower value function $W(t,x)$\ is of at most linear growth it belongs to $\Theta$,
   and so $W(t,x)$ is the unique viscosity solution in $\Theta$ of equation (6.1). Similarly we get that the upper value function
   $U(t,x)$\ is the unique viscosity solution in $\Theta$ of equation (5.2).\er

\br\mbox{  } If the Isaacs' condition holds, that is, if for all
$(t, x, y, p, X)\in [0, T]\times {\mathbb{R}}^n \times
{\mathbb{R}}\times {\mathbb{R}}^n\times {\mathbf{S}}^n ,$
$$H^-(t, x, y, p, X)=H^+(t, x, y, p, X),$$
then the equations (6.1) and (5.2) coincide and from the
uniqueness in $\Theta$ of viscosity solution it follows that the
lower value function $W(t,x)$ equals to the upper value function
$U(t,x)$ which means the associated stochastic differential game
has a value.\er

 \br\mbox{  } Let us assume that the coefficient of BSDE (4.5) $f(t,x,y,z,u,v)\equiv f(t,x,u,v)$\ is
 independent of $(y, z),$\ and denote by $\tilde{W}(t,x)$\ (resp.,
 $\tilde{U}(t,x)$)\ the lower value function (resp., the upper value function) defined by Fleming and Souganidis~\cite{FS1}, see Remark 4.2.
It is shown in~\cite{FS1} that $\tilde{W}(t,x)$\ is a viscosity
solution in $\Theta$ of (6.1) and $\tilde{U}(t,x)$\ a viscosity
solution in $\Theta$ of (5.2). Then, due to Theorem 6.1,
$W(t,x)=\tilde{W}(t,x)$\ and $U(t,x)= \tilde{U}(t,x), \ (t, x)\in
[0, T]\times {\mathbb{R}}^n.$\ Moreover, if the Isaacs' condition
holds then $W(t,x)=\tilde{W}(t,x)= \tilde{U}(t,x)=U(t,x).$\er

\end{document}